\documentclass[final,3p,times]{elsarticle}

\usepackage{amssymb}

\usepackage{lipsum}
\usepackage{amsfonts}
\usepackage{graphicx}
\usepackage{epstopdf}
\usepackage{algorithmic}
\ifpdf
  \DeclareGraphicsExtensions{.eps,.pdf,.png,.jpg}
\else
  \DeclareGraphicsExtensions{.eps}
\fi

\usepackage{natbib}
\usepackage{graphicx}
\usepackage{float}
\usepackage{subfig}
\usepackage{amsopn}
\usepackage{amsmath}
\usepackage{hyperref}
\usepackage{cleveref}
\usepackage{color,xcolor}
\usepackage{ulem}

\newcommand{\tpr}{\boldsymbol{\otimes}}   
\newcommand{\dpr}{\boldsymbol{\cdot}}     

\newcommand{\vc}[1]{\boldsymbol{#1}}      

\newcommand{\dudx}[2]{\frac{\partial{#1}}{\partial{#2}}}

\newtheorem{definition}{Definition}
\newtheorem{proposition}{Proposition}
\newtheorem{proof}{Proof}

\begin{document}

\begin{frontmatter}



\title{A reduced model for compressible viscous heat-conducting multicomponent flows}

\author[label1]{Chao Zhang}
\author[label1,label5]{Lifeng Wang}
\author[label1,label5]{Zhijun Shen}
\author[label1]{Zhiyuan Li}
\author[label3,label4]{Igor Menshov}

\address[label1]{Institute of Applied Physics and Computational Mathematics, Beijing, China}
\address[label3]{Keldysh Institute for Applied Mathematics RAS, Moscow, Russia}
\address[label4]{SRISA RAS, Moscow, Russia}
\address[label5]{Center for Applied Physics and Technology, HEDPS, Peking University, Beijing, China}

\begin{abstract}
In the present paper we propose a reduced temperature non-equilibrium model for simulating multicomponent flows with inter-phase heat transfer, diffusion processes (including the viscosity and the heat conduction) and external energy sources. We derive three equivalent formulations for the proposed model. 
The first formulation consists of balance equations for partial densities, the mixture momentum, the mixture total energy, and phase volume fractions. The second formulation is symmetric and obtained by replacing the equations for the mixture total energy and volume fractions in the first formulation with balance equations for the phase total energy. Replacing one of the phase total energy equation of the second formulation with the mixture total energy equation gives the third formulation. All the three formulations assume velocity and pressure equilibrium across the material interface. These equivalent forms provide different physical perspectives and numerical conveniences. Temperature equilibration and continuity across the material interfaces are achieved with the instantaneous thermal relaxation. Temperature equilibrium is maintained during the heat conduction process. The proposed models are proved to respect the thermodynamical laws. For numerical solution, the model is split into a hyperbolic partial differential equation (PDE) system and parabolic PDE systems. The former is solved with the high-order Godunov finite volume method that ensures the pressure-velocity-temperature (PVT) equilibrium conduction.  The parabolic PDEs are solved with both the implicit and the explicit locally iterative method (LIM) based on Chebyshev parameters. Numerical results are presented for several multicomponent flow problems with diffusion processes. Furthermore, we apply the proposed model to simulate the target ablation problem that is of significance to inertial confinement fusion. Comparisons with  one-temperature models in literature demonstrate the ability to maintain the PVT property and  superior convergence performance of the proposed model in solving multicomponent problems with diffusions. 
\end{abstract}



\begin{keyword}


Compressible multiphase flow  \sep viscosity \sep heat conduction \sep Godunov method \sep Chebyshev locally iterative method 
\end{keyword}

\end{frontmatter}


\section{Introduction}
\label{sec:intro}

The present research is devoted to the numerical modeling of compressible multicomponent flows including viscosity, heat conduction and external energy sources. This topic is of significance for applications in many fields, such as the inertial confinement fusion (ICF), the astrophysical events, detonation physics, oil/gas spilling, and so forth. 

Our work is performed in the framework of the diffuse interface method (DIM). The DIM captures the material interface by allowing artificial mixture of the fluids that is described in a thermodynamically consistent way. When simulating multiphase flows with resolved interfaces, one important property for DIMs is the capability to correctly simulate the pure translation of the isolated contact discontinuity between fluids with different equations of state (EOS). If temperature, pressure and velocity are uniformly distributed throughout the computational domain, their profiles should be maintained during the interface translation. This property is termed as the PVT property and formulated as follows:

\begin{definition}
An interface-capturing numerical scheme has the PVT property if it ensures  \[u_{i}^{n+1} = u = \text{const}, \; p_{i}^{n+1} = p = \text{const}, \; T_{i}^{n+1} = T = \text{const},\]
providing that
\[u_{i}^{n} = u = \text{const}, \; p_{i}^{n} = p = \text{const}, \; T_{i}^{n} = T = \text{const},\]
\end{definition}
here, $u$, $p$ and $T$ are the mixture velocity, the pressure and the temperature, respectively. The subscript $i$ and superscript $n$ denote the spatial index of the computational cell and solution time step, respectively. Spurious oscillations in pressure and temperature arise once this PVT property is violated \cite{abgrall1996prevent,Thornber2018,Johnsen2012Preventing}.

In literature, the most widely used models for implementing DIM include the two-phase Baer-Nuziato (BN) model \cite{BAER1986861}, its variant \cite{saurel1999multiphase} or  its reduced models \cite{kapila2001two,kreeft2010new}.
The BN model is originally developed for the deflagration-to-detonation transition (DDT) problem. This model is then extended to and widely used in the simulation of compressible multiphase/multicomponent flows. 
In this model each phase is described with a complete set of their own flow variables (density, velocity, pressure and temperature). The phase interactions happen within the diffused zone. Outside this diffused zone, the flow behaviours of each pure component are governed by Euler equations (or Navier-Stokes equations with viscous terms being included). The phase interaction is described by relaxation terms driving the flow variables towards an equilibrium state.

The BN model ensures the PVT property. It is  hyperbolic, physically complete and thermodynamically consistent. However, its numerical solution is quite complicated due to its complex wave structure and stiff relaxations. For simplicity, various reduced models are derived, for example, the four-equation model \cite{Lemartelot2014}, the five-equation model \cite{kapila2001two}, and the six-equation model \cite{saurel2009simple}. Among these reduced models, the one-temperature four-equation model is the most straightforward to consider the heat conduction terms and widely used for simulating boiling and combustion problems where the inter-phase heat exchange takes place violently. However, it is well known that this model violates the PVT property and results in spurious oscillations in pressure and erroneous spikes in temperature in the vicinity of material interfaces \cite{Johnsen2012Preventing}. The erroneous temperature spike is then spread out throughout the computational domain by the heat conduction terms, resulting in convergence problems.  This problem arises as a result of incompatibility of the isothermal closure with the interface jump conditions, i.e., continuity in pressure and normal velocity.

Kapila et al. \cite{kapila2001two} have derived a five-equation model for two-phase flows without diffusions by performing asymptotic analysis of the BN model in the limit of instantaneous mechanical (velocity and pressure) relaxations. In this model, governing equations are written for the partial densities, the mixture momentum, the mixture total energy and the volume fraction. This model is free of the pressure oscillation problems when computing interface problems \cite{murrone2005five}, yet much simpler than the complete BN model. It allows non-equilibrium phase temperatures and two phase entropies.  It is thermodynamically consistent and ensures the interface jump conditions.  Moreover, this reduced model captures the dynamically arising material interfaces, for example, the cavitation interfaces. 

In view of the above discussions, we are more interested in the temperature non-equilibrium model. To the authors' knowledge, the work including diffusions (viscosity and heat conduction) into the temperature non-equilibrium five-equation model is absent in literature so far. The main difficulty consists in how to describe the inter-phase energy exchange with only one energy equation being included. In fact, in Kapila's five-equation model \cite{kapila2001two} the temperature relaxation between phases is described implicitly by including a right-hand side (RHS) term related to the temperature relaxation in the volume fraction equation. However, different from the temperature relaxation terms (that are eliminated after summing), the heat conduction terms appear in the mixture energy equation. As a result, we obtain one energy equation with two phase temperatures, whose solution is not straightforward, especially when finite temperature relaxation rate is considered. Therefore,  to solve this problem, we seek a different model formulation.

Instead of the original Kapila's formulation (with one mixture energy equation) of the five-equation model, an equivalent formulation consisting of two phase total energy equations (without diffusion processes) has been derived on the basis of physical laws in \cite{kreeft2010new}.
Inspired by this formulation, we derive a new formulation consisting of multiple phase energy equations and including the diffusion processes. Our derivation procedure is based on thermodynamical relations, which is different from that of \cite{kreeft2010new}. To ensure thermodynamical consistency, we first write a BN-type model for $N$-phase flows with diffusions, then perform asymptotic analysis in the limit of instantaneous mechanical relaxations. The asymptotic analysis leads to three equivalent formulations of the reduced model. The first formulation consists of balance equations for partial densities, the mixture momentum, the mixture total energy, and volume fractions. The second formulation is symmetric and obtained by replacing the equations for the mixture total energy and volume fractions in the first formulation with balance equations for the phase total energy  with interaction and energy-exchange terms. Replacing one of the phase total energy equation of the second formulation with the mixture total energy equation gives the third form. With these model formulations, we can consider the energy exchange explicitly. We prove that the proposed model ensures the PVT property and the entropy inequality.

We develop numerical methods for solving the proposed model with the fractional step method. In the numerical implementation, the model is split into four sub-systems, i.e., the hydrodynamic part, the viscous part, the temperature relaxation part, and the heat conduction part. For the hydrodynamic part, we use the single-energy equation. For the rest parts, multiple phase total energy equations are calculated. The homogeneous hyperbolic equation (i.e., the hydrodynamic part) is   
solved with the Godunov finite volume method. The diffusion processes (including viscous and heat conduction terms) are governed by parabolic PDEs, that are solved with the locally iterative method based on the Chebyshev parameters.
The heat conduction equations are solved maintaining the temperature equilibrium, i.e., assuming an instantaneous temperature relaxation.  Note that the finite temperature relaxation can also be considered straightforwardly with the proposed new formulations.

The rest of the present paper is organized as follows. In \Cref{sec:model} we derive three equivalent reduced models by performing the asymptotic analysis on the BN-type model in the limit of instantaneous mechanical relaxations. In \Cref{sec:numer_meth} we develop numerical methods for solving the proposed model. In \Cref{sec:numer_res} we present numerical  results for several multicomponent problems with diffusions and apply the model and numerical methods to the target ablation  problem in the field of ICF.

\section{Model formulation}
\label{sec:model}
\subsection{The parent model with viscosity, heat conduction and external heat source}
On the basis of the BN model \cite{BAER1986861} or its variant for compressible multiphase flows \cite{Romenski2007Conservative,saurel1999multiphase,saurel2001multiphase}, the three-phase model \cite{HERARD2007732}, and the basic continuum mechanics for multiphase flows in \cite{Nigmatulin1987}, we give a general form for the $N$-phase model  as follows:
\begin{subequations} \label{eq:bn}
\begin{align}
\label{eq:bn:mass}
\dudx{\alpha_k\rho_k}{t} + \nabla\dpr(\alpha_k\rho_k\vc{u}_k) = 0, \\
  \label{eq:bn:mom}
  \dudx{\alpha_k\rho_k\vc{u}_k}{t} + \nabla\dpr\left(\alpha_k\rho_k\vc{u}_k\tpr\vc{u}_k - \alpha_k \overline{\overline{T}}_k \right) = \sum_{j\neq k}^N {\mathcal{P}}_{jk},\\
  \label{eq:bn:en}
  \dudx{\alpha_k \rho_k E_k}{t} + 
  \nabla\dpr\left(
\alpha_k \rho_k E_k \vc{u}_k   - \alpha_k \overline{\overline{T}}_k \dpr \vc{u}_k
  \right) 
  = 
 \sum_{j\neq k}^N \mathcal{E}_{jk} + q_k + \mathcal{I}_k,\\
   \label{eq:bn:vol}
  \dudx{\alpha_l }{t} + 
  \vc{u}_I \dpr \nabla\alpha_l = \sum_{j\neq l}^{N}\mathcal{F}_{jl},
\end{align}
\end{subequations}
where the notations used are standard: $\alpha_k, \; \rho_k, \; \vc{u}_k, \; p_k, \; \overline{\overline{T}}_k, \; E_k$ are the volume fraction, phase density, velocity, pressure, stress tensor, and total energy of component $k$. Due to the saturation constraint for volume fractions 
$\sum \alpha_k  = 1$, the last equations are written only for only $N-1$ volume fractions, where $N$ is the number of components. We define  $k\in \Psi = \{1,2,\cdots,N\},\; l\in \Psi^{\prime} = \{2,\cdots,N\} \subset \Psi.$  The phase density $\rho_k$ is defined as the mass per unit volume occupied by the phase $k$. The mixture density $\rho$ is the sum of the partial densities $\alpha_k \rho_k$, i.e., $\rho = \sum {\alpha_k \rho_k}$. 
 The total energy is $E_k  = e_k + \mathcal{K}_k$ where $e_k$, and $\mathcal{K}_k = \frac{1}{2}\vc{u}_k \dpr \vc{u}_k$ are the  internal energy and kinetic energy, respectively.

The phase stress tensor, $\overline{\overline{T}}_k$, can be written as
\begin{equation}
 \overline{\overline{T}}_k  = - p_k \overline{\overline{I}} +  \overline{\overline{\tau}}_k.
\end{equation}

For the viscous part we use the  Newtonian approximation
\begin{equation}\label{eq:newton_vis}
\overline{\overline{\tau}}_k = 2\mu_k \overline{\overline{D}}_k + \left(\mu_{b,k} - \frac{2}{3}\mu_k \right) \nabla \dpr \vc{u}_k,
\end{equation}
where $\mu_k > 0$ is the coefficient of shear viscosity and  $\mu_{b,k} > 0$ is the coefficient of bulk viscosity, $\overline{\overline{D}}_k$ is defined as
\[
\overline{\overline{D}}_k = \frac{1}{2} \left[ \nabla \vc{u}_k + \left(  \nabla \vc{u}_k \right)^{\text{T}} \right].
\]

The terms $\mathcal{P}_{jk}$ and $\mathcal{E}_{jk}$ are the momentum and energy exchange intensity between the phase $j$ and the phase $k$, respectively,
\begin{equation}
\mathcal{P}_{jk}  =  \overline{\overline{T}}_{jk,I} \dpr \nabla {\alpha_j}  + \mathcal{M}_{jk},
\end{equation}
\begin{equation}
\mathcal{E}_{jk} = \vc{u}_{I} \dpr \left( \overline{\overline{T}}_{jk,I} \dpr \nabla {\alpha_j} \right) 
    + \vc{u}_I \dpr \mathcal{M}_{jk} 
- p_{jk,I} \sum_{l\neq j} \mathcal{F}_{lj} + \mathcal{Q}_{jk}.
\end{equation}

The inter-phase relaxation terms include the velocity relaxation $\mathcal{M}_{jk}$, the pressure relaxation $\mathcal{F}_{jk}$, and the temperature relaxation $\mathcal{Q}_{jk}$. They are as follows:
\begin{equation}\label{eq:relaxations}
\begin{split}
\mathcal{M}_{jk}= \vartheta_{jk} \left( \vc{u}_{j}-\vc{u}_k \right), \quad \mathcal{F}_{jk}=\eta_{jk} \left( {p}_{k} - {p}_{j} \right), \quad \mathcal{Q}_{jk}=\varsigma_{jk} \left( T_{j} - T_{k} \right).
\end{split}
\end{equation}
where $\vartheta_{jk} = \vartheta_{kj} > 0$, $\eta_{jk} = \eta_{kj} > 0$, and $\varsigma_{jk} = \varsigma_{kj} > 0$ are the corresponding relaxation rates between  phase $j$ and phase $l$.

The variables with the subscript ``I'' represent the variables at interfaces, for which there are several possible definitions \cite{saurel1999multiphase,perigaud2005compressible,saurel2018diffuse}. Whatever definitions we choose, $\lim_{\eta\to\infty} p_{jk,I} = \lim_{\eta\to\infty}p_{kj,I} = p$, $\lim_{\vartheta\to\infty}\vc{u}_{I} = \lim_{\vartheta\to\infty}\vc{u}_{k} = \vc{u}$, and $\lim_{\vartheta\to\infty} {\overline{\overline{\tau}}}_{jk,I} = \lim_{\vartheta\to\infty} {\overline{\overline{\tau}}}_{kj,I} = \lim_{\vartheta\to\infty}{\overline{\overline{\tau}}}_{k} = {\overline{\overline{\tau}}}$. The stress tensor on the interface between the phase $j$ and the phase $k$  is
$\overline{\overline{T}}_{jk,I} = - p_{jk,I} \overline{\overline{I}} + \overline{\overline{\tau}}_{jk,I}$ .

Moreover, we have
\begin{equation}
\mathcal{P}_{jk} = \mathcal{P}_{kj}, \; \mathcal{P}_{kk} = 0, \;
\mathcal{E}_{jk} = \mathcal{E}_{kj}, \; \mathcal{E}_{kk} = 0.
\end{equation}

Here, we introduce different interfacial pressures for each pair of phases, $p_{jk,I}$. However,  we only define one interfacial velocity. As the analysis in \cite{HERARD2007732} shows, defining multiple interfacial velocities may result in the violation    of the maximum principle for the volume fractions. The interfacial parameters should be defined in such a way that the hyperbolicity and thermodynamical consistency of the model are ensured.

The heat conduction term is given as
\begin{equation}\label{eq:q}
q_k = \nabla \dpr \vc{q}_k,
\end{equation}
where $\vc{q}_k = - \alpha_k \lambda_k \nabla T_k$ is the Fourier's heat flux.

The external heat source is
\begin{equation}\label{eq:I}
\mathcal{I}_k = \alpha_k I_k,
\end{equation}
where $I_k$ denotes the the intensity of the external heat source released in the phase $k$, $I_k \geq 0$.

Each phase obeys their own stiffened gas (SG) EOS, which takes the following form:
\begin{equation}\label{eq:eos}
\rho_{k} e_{k}=\frac{p_{k}+\gamma_{k} p_{\infty, k}}{\gamma_{k}-1}+\rho_{k} w_{k} = \rho_{k} C_{v, k} T_{k}+p_{\infty, k}+\rho_{k} w_{k},
\end{equation}
where  $C_{v,k}$ is the specific heat at constant volume. Parameters $\gamma_k$, $p_{\infty,k}$, and $w_k$ are constants depending on the properties of the component $k$.

In the two-phase flow ($N=2$), if the interface pressure and velocity are chosen to be those of the more compressible and less compressible phases, respectively, the model (\ref{eq:bn}) reproduces the BN seven-equation model \cite{BAER1986861}. Without the diffusion processes, the seven-equation model is  hyperbolic with the following set of wave speeds $u_k \pm a_k, u_k, u_I$, where $a_k$ is the sound speed 
\begin{equation}
a_k ^2 = \left( \dudx{p_k}{\rho_k} \right)_{s_k} = \frac{\frac{p_k}{\rho_k^2} - \left( \dudx{e_k}{\rho_k} \right)_{p_k} }{\left( \dudx{e_k}{p_k} \right)_{\rho_k}} > 0.
\end{equation}


The three-phase model in \cite{HERARD2007732} is also a particular case of (\ref{eq:bn}) when $N=3$ and in the absence of the diffusion terms. In this model the interface pressures are defined as
\begin{equation}
\begin{array}{l}
p_{11,I}=0, \quad p_{12,I}=p_{2}, \quad p_{13,I}=p_{3}, \\
p_{21,I}=p_{2}, \quad p_{22,I}=0, \quad p_{23,I}=p_{2} \\
p_{31,I}=p_{3}, \quad p_{32,I}=p_{3}, \quad p_{33,I}=0 .
\end{array}
\end{equation}

This definition of interface pressures ensures the hyperbolicity of the equation system (\ref{eq:bn}) and is consistent with the entropy inequality, as demonstrated in \cite{HERARD2007732}.

\subsection{Equations for the primitive variables}
In this section, we derive equations for some primitive variables, which are to be used for further analysis. We introduce the material derivative related to the phase velocity $\mathbf{u}_k$ and the interfacial velocity $\mathbf{u}_I$,
\begin{equation}
 \frac{\text{D}_g \Phi}{\text{D} \Phi} = \dudx{\Phi}{t} + \mathbf{u}_g \cdot \nabla{\Phi}, \;\; g = k, I.
\end{equation} 

Some useful relations are given as follows:
\begin{equation}\label{eq:gibbs}
T_k \frac{\text{D}_k s_k}{\text{D} t} = \frac{\text{D}_k e_k}{\text{D} t} - \frac{p_k}{\rho_k^2} \frac{\text{D}_k \rho_k}{\text{D} t}, \;\;\frac{\text{D}_k e_k}{\text{D} t} = \chi_k \frac{\text{D}_k \rho_k}{\text{D} t} + \xi_k \frac{\text{D}_k p_k}{\text{D} t}, \;\;
\frac{\text{D}_k p_k}{\text{D} t} = a_k^2 \frac{\text{D}_k \rho_k}{\text{D} t} + \omega_k \frac{\text{D}_k s_k}{\text{D} t},
\end{equation}
where the first expression is the Gibbs relation, $s_k$ is the phase entropy. The second and third equations are total differentials of $e_k$ and $p_k$, 
\[\chi_k = \dudx{e_k}{\rho_k}\Big|_{p_k}, \; \xi_k = \dudx{e_k}{p_k}\Big|_{\rho_k}, \;\omega_k = \dudx{p_k}{s_k}\Big|_{\rho_k},\]
simple manipulations of \cref{eq:gibbs} lead to  $\chi_k = p_k/\rho_k^2 - \xi_k a_k^2$.

Following a procedure similar to \cite{murrone2005five,zhangPHD2019}, one can obtain the following equations with respect to the primitive variables (including the phase entropy $s_k$, the phase velocity $\mathbf{u}_k$, the phase pressure $p_k$ and volume fraction $\alpha_k$):
\begin{subequations} \label{eq:bn_prim}
\begin{align}
\alpha_{k} \rho_{k} T_{k} \frac{\mathrm{D}_{k} s_{k}}{\mathrm{D} t} = \left(\boldsymbol{u}_{I}-\boldsymbol{u}_{k}\right) \cdot \left[ \mathcal{M}_{k}^{sum}  - \sum_{j\neq k}^N \left(p_{jk,I}-p_{k}\right) \cdot \nabla \alpha_{j}  +  \sum_{j\neq k}^{N} \left(  \overline{\overline{\tau}}_{jk,I} \cdot \nabla \alpha_j \right) \right] + \sum_{j\neq k}^N \left(p_{jk,I}-p_{k}\right) {\mathcal{F}}_{j}^{sum}
+\mathcal{G}_{k} \label{eq:sk}   \\
\alpha_{k} \rho_{k} \frac{\mathrm{D}_{k} \boldsymbol{u}_{k}}{\mathrm{D} t} = \nabla\cdot \left( \alpha_k \overline{\overline{T}}_k \right)  +  {\mathcal{P}_{k}^{sum}} \label{eq:uk}\\
\frac{\mathrm{D}_{k} p_{k}}{\mathrm{D} t}=  \sum_{j\neq k}^{N} \frac{\rho_{k} a_{jk,I}^{2}}{\alpha_{k}} {\mathcal{F}}_{j}^{sum} + \frac{\mathbf{u}_I - \mathbf{u}_k}{\alpha_k \rho_k \xi_k} \left[ \sum_{j\neq k}^{N} \left( \overline{\overline{\tau}}_{jk,I} - \xi_k \rho_k^2 a_{jk,I}^2 \overline{\overline{I}} \right) \cdot \nabla\alpha_j + \mathcal{M}_k^{sum} \right]  + \frac{\Gamma_k \mathcal{G}_k }{ \alpha_k } - A_k \nabla\cdot \mathbf{u}_k \label{eq:pk}\\
\frac{\mathrm{D}_{I} \alpha_{k}}{\mathrm{D} t}={\mathcal{F}}_{k}^{sum} \label{eq:alpk}
\end{align}
\end{subequations}
where \[\frac{\rho_{k} a_{jk,I}^{2}}{\alpha_{k}} = \frac{\rho_{k} a_{k}^{2}}{\alpha_{k}} + \frac{p_{jk,I}-p_{k}}{\alpha_{k} \rho_{k} \xi_{k}}, \;\; \Gamma_k = \frac{1}{\rho_k \xi_k}, \;\; A_k = \rho_k a_k^2,\]
\begin{equation}\label{eq:Gk}
\mathcal{G}_k = \alpha_k \overline{\overline{\tau}}_k : \overline{\overline{D}}_k + \mathcal{Q}_k^{sum} + q_k + \mathcal{I}_k,
\end{equation}
\[
{\Phi}_{k}^{sum} = \sum_{j\neq k}^N {\Phi}_{jk}, \;\; {\Phi} = \mathcal{M}, \; \mathcal{F}, \; \mathcal{Q}, \; \mathcal{P}. 
\]

One can check that in the case $N=2$ and in the absence of the diffusion and external energy terms, \cref{eq:sk,eq:uk,eq:pk,eq:alpk} are reduced to those in \cite{murrone2005five}.

\subsection{Model reduction}
Although physically complete, due to the complex wave structure and stiff relaxations of the BN model, its numerical implementation is rather cumbersome. According to the evaluations of different mechanical/thermal relaxation time scales in \cite{kapila2001two}, this model can be reduced significantly for certain application scenarios. A variety of reduced models consisting of three to six equations are proposed in literature \cite{lund2012hierarchy,kapila2001two,Lemartelot2014,murrone2005five,lund2012hierarchy}. 

One of the best-known reduced models is the Kapila's five-equation model \cite{kapila2001two} that is developed for modelling the deflagration-to-detonation transition in granular materials. This model is developed for two-phase flows in the absence of viscosity, heat conduction and external source terms. It is a limit model of the BN model in the case of instantaneous mechanical (pressure and velocity) relaxations. This model is then widely used for simulating compressible flows with resolved fluid interfaces thanks to its simplicity and ability to ensure the interface jump conditions (i.e., the continuity in pressure and normal velocity).

In the present work, we are interested in multicomponent flows with resolved material interfaces, therefore, the instantaneous mechanical relaxation assumption is also adopted here, i.e.,
\begin{equation}
\vartheta_{jk} \to \infty, \; \eta_{jk} \to \infty,
\end{equation}
or
\begin{equation}
\epsilon = 1 / \vartheta_{jk} \approx 1 / \eta_{jk}  \to 0,
\end{equation}
where $\epsilon$ denotes the mechanical relaxation time.

\subsubsection{Velocity equilibration}
We assume the following asymptotic expansion for velocity
\begin{equation}\label{eq:asymp_u}
\vc{u}_k = \vc{u}_k^{(0)} + \epsilon \vc{u}_k^{(1)} + \mathcal{O}(\epsilon^2).
\end{equation}

Inserting \cref{eq:asymp_u} into \cref{eq:uk} yields $\vc{u}_1^{(0)} = \vc{u}_2^{(0)} = \cdots = \vc{u}_N^{(0)} = \vc{u}^{(0)}$ in the order  $\mathcal{O}(\vartheta_{jk}) = \mathcal{O}(1/\epsilon)$. This means that the velocities are in equilibrium in $\mathcal{O}(1)$ order, which significantly simplifies the wave structure of the model.

\subsubsection{Pressure equilibration}
For pressure, we similarly assume the asymptotic expansion 
\begin{equation}\label{eq:asymp_p}
p_k = p_k^{(0)} + \epsilon p_k^{(1)} + \mathcal{O}(\epsilon^2).
\end{equation}

Inserting \cref{eq:asymp_p} into \cref{eq:pk}, in the order of $\mathcal{O}(1/\epsilon)$ we obtain
\begin{equation}\label{eq:pres_eq_0}
p_1^{(0)} =p_2^{(0)} = \cdots =p_N^{(0)} = p^{(0)}.
\end{equation}

In the order of $\mathcal{O}(1)$, the phase pressure equation \cref{eq:pk} is reduced to
\begin{equation}\label{eq:pk_o1}
\frac{\mathrm{D}_{k} p_{k}^{(0)}}{\mathrm{D} t} = \sum_{j\neq k}^{N} \frac{A_k}{\alpha_k} \mathcal{F}_j^{sum,(0)} + \frac{\Gamma_k \mathcal{G}_k}{\alpha_k} - A_k \nabla\dpr \vc{u}^{(0)}.
\end{equation}

Combining \cref{eq:pres_eq_0,eq:pk_o1} and 
\begin{equation}
\sum_{j=1}^{N} \mathcal{F}_j^{sum,(0)} = 0,
\end{equation}
one can obtain the following solution for $\mathcal{F}_k^{sum,(0)}$:
\begin{equation}\label{eq:FK0}
\mathcal{F}_k^{sum,(0)} = \alpha_{k} \frac{A-A_{k}}{A_{k}} \nabla\dpr \vc{u}^{(0)} + A \sum_{j \neq k}^{N} \frac{{\mathcal{G}}_{k} \alpha_{j}-{\mathcal{G}}_{j} \alpha_{k}}{A_{k} A_{j}},
\end{equation}
where \[\frac{1}{A} = \sum_{k=1}^{N} \frac{\alpha_k}{A_k}.\]

Retaining only to the order $\mathcal{O}\left( 1\right)$, the volume fraction equation becomes
\begin{equation}
\dudx{\alpha_k}{t} + \vc{u}^{(0)} \dpr \nabla \alpha_k  =  \mathcal{F}_k^{sum,(0)}.
\end{equation}

\subsubsection{The final model}
We summarize the derived one-velocity one-pressure model in the order of $\mathcal{O}(1)$ as follows:
\begin{subequations} \label{eq:reduced_five}
\begin{align}
\dudx{\alpha_k \rho_k}{t} + \nabla \dpr \left( {\alpha_k \rho_k {\vc{u}}} \right) = 0,  \label{eq:reduced_five:mass} \\
\dudx{\rho{\vc{u}}}{t} + \nabla\dpr\left(\rho{\vc{u}}\tpr{\vc{u}} - \overline{\overline{T}} \right)  = 0, \label{eq:reduced_five:mom} \\
\dudx{ \rho E}{t} + 
  \nabla\dpr\left(
\rho E {\vc{u}}   - \overline{\overline{T}} \dpr {\vc{u}}
  \right) 
  =  \sum_{k=1}^{N} \left( q_k + {\mathcal{I}}_k \right), \label{eq:reduced_five:en} \\
\dudx{\alpha_l}{t} + {\vc{u}} \dpr \nabla \alpha_l  =  \mathcal{F}_l^{sum,(0)}, \label{eq:reduced_five:vol}
\end{align}
\end{subequations}
where the superscript ``${(0)}$'' over the velocity and pressure is omitted for simplicity. The equation for the mixture momentum (\ref{eq:reduced_five:mom}) and that for the mixture total energy (\ref{eq:reduced_five:en}) are obtained by summing 
\cref{eq:bn:mom} and \cref{eq:bn:en}, respectively. The volume fraction equation (\ref{eq:reduced_five:vol}) is written for $N-1$ components, the remaining one volume fraction can be solved with the saturation condition for volume fractions. The model consists of $2N+1$ equations in the case of $N$ phases. Note that in the case $N=2$ and in the absence of diffusions and external energy terms, this model is reduced to the five-equation one \cite{murrone2005five,kapila2001two}.

Regarding the viscous terms, our model is different from that of Perigaud et al. \cite{perigaud2005compressible}. The difference consists in that in our model the dissipated part of mechanical energy  (i.e., $\alpha_k \overline{\overline{\tau}}_k : \overline{\overline{D}}_k$) has an impact on the volume fraction evolution, which can be seen from \cref{eq:reduced_five:vol,eq:FK0,eq:Gk}. This term is absent in \cite{perigaud2005compressible}.

\subsection{New formulations of the reduced model}
In this section, we derive new formulations that are equivalent to  the model (\ref{eq:reduced_five}). These new formulations explicitly demonstrate the interaction and energy exchange terms between different phases, making them very straightforward to consider the thermal relaxation and the heat conduction processes. 

Similar to the formulation in Kreeft et al. \cite{kreeft2010new}, our formulation includes the phase energy equations with inter-phase exchange terms and excludes the volume fraction equations. Our derivation is based on the reformulation of \cref{eq:reduced_five} with thermodynamic relations.

In the limit of zero mechanical relaxation time, the evolution equation for the phase entropy (\ref{eq:sk}) is reduced to 
\begin{equation}
\alpha_k \rho_k T_k \frac{\text{D} s_k}{\text{D} t}= \mathcal{G}_k.
\end{equation}

Invoking the Gibbs relation
\begin{equation}
\mathrm{d} e_{k} = T_{k} \mathrm{~d} s_{k} + \frac{p}{\rho_{k}^{2}} \mathrm{~d} \rho_{k},
\end{equation}
one can deduce
\begin{equation}\label{eq:De}
\alpha_k \rho_k \frac{\text{D} e_k}{\text{D} t} = \mathcal{G}_k + \frac{\alpha_k p}{\rho_{k}} \frac{\text{D} \rho_k}{\text{D} t}.
\end{equation}

Combination of \cref{eq:reduced_five:mass,eq:reduced_five:vol} leads to
\begin{equation}\label{eq:Drho}
\frac{\text{D} \rho_k}{\text{D} t} = - \frac{\rho_k}{\alpha_k} \mathcal{F}_k^{sum,(0)} - \rho_k \nabla\dpr \mathbf{u}.
\end{equation}

By using \cref{eq:De,eq:Drho,eq:reduced_five:mom}, one can obtain
\begin{equation}\label{eq:Ek}
\dudx{\alpha_k \rho_k E_k}{t} + \nabla \dpr \left( \alpha_k \rho_k E_k \mathbf{u} + p\mathbf{u}\right) = \mathcal{G}_{D,k} + \mathcal{G}_{M,k} + \mathcal{G}_{T,k} + \mathcal{I}_k + \mathcal{Q}_{k}^{sum} ,
\end{equation}
where 
\begin{subequations}
\begin{align}
\mathcal{G}_{D,k} &= \nabla \dpr \left( \alpha_k \overline{\overline{\tau}}_k \dpr \mathbf{u}  \right) +q_{k}  , \\
\mathcal{G}_{M,k} &= \left( \alpha_k - y_k \right) \mathbf{u} \dpr \nabla p + p\mathbf{u} \dpr \nabla \alpha_k, \quad y_k = \alpha_k \rho_k / \rho \\ 
\mathcal{G}_{T,k} &= - p  \mathcal{F}_k^{sum,(0)}.
\end{align}
\end{subequations}

The terms $\mathcal{G}_{D,k}$, $\mathcal{G}_{M,k}$ and $\mathcal{G}_{T,k}$ have clear physical explanations:
\begin{enumerate}
\item[•] The term $\mathcal{G}_{D,k}$ represents the phase total energy change rate due to diffusion processes: the viscous diffusion $\nabla \dpr \left( \alpha_k \overline{\overline{\tau}}_k \dpr \mathbf{u}  \right)$, and the heat conduction $q_k$. 
\item[•] The term $\mathcal{G}_{M,k}$ is the mechanical work rate. The first part $\left( \alpha_k - y_k \right) \mathbf{u} \dpr \nabla p$  represents the the work rate of the force, exerted on the fluid $k$ to maintain velocity equilibrium. The second part $p\mathbf{u} \dpr \nabla \alpha_k$  represents the work rate of the pressure force acting on the material interface due to the spatial variation of the flow area (volume fraction).
\item[•] The term $\mathcal{G}_{T,k}$ is the rate of thermodynamic work. This term maintains the pressure equilibrium in the expansion or compression of the phase $k$. 
\end{enumerate}

With \cref{eq:Ek} we can present two formulations equivalent to  \cref{eq:reduced_five}:
\begin{enumerate}
\item[(1)] The symmetric formulation
\begin{subequations} \label{eq:reduced_five_1}
\begin{align}
\dudx{\alpha_k \rho_k}{t} + \nabla \dpr \left( {\alpha_k \rho_k \vc{u}} \right) = 0,  \label{eq:reduced_five:mass_1} \\
\dudx{\rho\vc{u}}{t} + \nabla\dpr\left(\rho\vc{u}\tpr\vc{u} + p \overline{\overline{I}} \right)  = 
\nabla \dpr \overline{\overline{\tau}}, \label{eq:reduced_five:mom_1} \\
\dudx{\alpha_k \rho_k E_k}{t} + \nabla \dpr \left( \alpha_k \rho_k E_k \mathbf{u} + p\mathbf{u} \right) =  \mathcal{G}_{D,k} + \mathcal{G}_{M,k} + \mathcal{G}_{T,k}  + \mathcal{I}_k + \mathcal{Q}_{k}^{sum}. \label{eq:reduced_five:enk_1} 
\end{align}
Here, \cref{eq:reduced_five:enk_1} is written for all the phases, $k=1,2,\cdots,N.$
\end{subequations}

\item[(2)] The non-symmetric formulation
\begin{subequations} \label{eq:reduced_five_2}
\begin{align}
\dudx{\alpha_k \rho_k}{t} + \nabla \dpr \left( {\alpha_k \rho_k \vc{u}} \right) = 0,  \label{eq:reduced_five:mass_2} \\
\dudx{\rho\vc{u}}{t} + \nabla\dpr\left(\rho\vc{u}\tpr\vc{u} + p \overline{\overline{I}} \right)  = 
\nabla \dpr \overline{\overline{\tau}}, \label{eq:reduced_five:mom_2} \\
 \dudx{\rho E}{t} +  \nabla\dpr \left( \rho E \vc{u} +  p \vc{u} \right) = \nabla \dpr \left( \overline{\overline{\tau}} \dpr \vc{u} \right) + \sum_{k=1}^{N} \left( q_k + \mathcal{I}_k \right), \label{eq:reduced_five:en_2} \\
\dudx{\alpha_l \rho_l E_l}{t} + \nabla \dpr \left( \alpha_l \rho_l E_l \mathbf{u} + p\mathbf{u} \right) = \mathcal{G}_{D,l} + \mathcal{G}_{M,l} + \mathcal{G}_{T,l}  + \mathcal{I}_l + \mathcal{Q}_{l}^{sum}. \label{eq:reduced_five:Ek2_2}
\end{align}
\end{subequations}
Here, \cref{eq:reduced_five:Ek2_2} is written for $N-1$ components.
\end{enumerate}

In the symmetric formulation (\ref{eq:reduced_five_1}), the mixture total energy equation (\ref{eq:reduced_five:en}) and the volume fraction equation (\ref{eq:reduced_five:vol}) of the original formulation are replaced by the phase total energy equations (\ref{eq:reduced_five:enk_1}). In the non-symmetric formulation we write the phase total energy equation (\ref{eq:reduced_five:Ek2_2}) instead of the volume fraction equation (\ref{eq:reduced_five:vol}). 

Although these two formulations  (\cref{eq:reduced_five_1,eq:reduced_five_2}) do not include the evolution equation for the volume fraction that describes the material interface topology, the volume fraction can be determined from the state variables $\left[ \alpha_k \rho_k, \; \rho \mathbf{u}, \; \alpha_k \rho_k E_k \right]$ or $\left[ \alpha_k \rho_k, \; \rho \mathbf{u}, \; \rho E, \; \alpha_l \rho_l E_l \right]$ by using the pressure equilibrium condition. Moreover, it is more straightforward to deal with the energy exchange terms with these new formulations in comparison with the original formulation with only one balance equation for the mixture total energy.

\subsection{Thermodynamical consistency}
We define the mixture entropy by assuming its mass additivity
\begin{equation}
\rho s = \sum_{k=1}^{N} m_k s_k, \quad \text{or} \quad s = \sum_{k=1}^{N} y_k s_k.
\end{equation}
where the partial density $m_k = \alpha_k \rho_k$.

Then by using the \cref{eq:sk,eq:reduced_five:mass}, we deduce
\begin{equation}
\rho \frac{\mathrm{D} s}{\mathrm{D} t} = \sum_{k=1}^{N}  \frac{\mathcal{G}_{k}}{T_{k}}.
\end{equation}

We then split the material derivative of the entropy into two parts:
\begin{equation}
\rho \frac{\mathrm{D} s}{\mathrm{D} t} = \rho \frac{\mathrm{D}^{(ext)} s}{\mathrm{D} t} + \rho \frac{\mathrm{D}^{(int)} s}{\mathrm{D} t},
\end{equation}
where the first term represents the entropy variation due to the entropy flux from external environment, and the second term is that due to the internal entropy production. The latter is non-negative according to the second law of thermodynamics. Here the entropy flux is due to the heat flux,
\begin{equation}\label{eq:ent_ext}
\rho \frac{\mathrm{D}^{(ext)} s}{\mathrm{D} t} = - \sum_{k=1}^{N} \nabla \dpr \left( \frac{\vc{q_k}}{T_k} \right),
\end{equation}
then the entropy production is 
\begin{equation}\label{eq:ent_int}
\rho \frac{\mathrm{D}^{(int)} s}{\mathrm{D} t} = \sum_{k=1}^{N}  \left( \frac{\alpha_{k} \overline{\overline{\tau}}_{k}: \overline{\overline{D}} + Q_{k}^{s u m} +\mathcal{I}_{k}}{T_k} - \frac{\vc{q}_k \dpr \nabla T_k}{T_k^2} \right).
\end{equation}

\begin{proposition}
The entropy production of the mixture entropy defined in the reduced model is non-negative, i.e.,
\begin{equation}\label{eq:ent_ineq}
\rho \frac{\mathrm{D}^{(int)} s}{\mathrm{D} t} \geq 0.
\end{equation}
\end{proposition}
\begin{proof}
We analyse the RHS terms \cref{eq:ent_int} one by one.

By simple tensor manipulations of \cref{eq:newton_vis}, one can prove that $\alpha_{k} \overline{\overline{\tau}}_{k}: \overline{\overline{D}} \geq 0$.

For the term including the temperature relaxation,
\begin{equation}
\sum_{k=1}^{N} \frac{Q_{k}^{s u m}}{T_k} = \sum_{j\neq l} \left(  \frac{\mathcal{Q}_{jl}}{T_l} +  \frac{\mathcal{Q}_{lj}}{T_j} \right) = \sum_{j\neq l} \frac{\varsigma_{lj} \left( T_l - T_j\right)^2}{T_l T_j} \geq 0.
\end{equation}
Moreover, $\mathcal{I}_k >0$,
\begin{equation}
\sum_{k=1}^{N} \frac{\mathcal{I}_k}{T_k} > 0 .
\end{equation}
Invoking the Fourier heat flux $\vc{q}_k = -\alpha_k \lambda_k \nabla T_k$, for the last term we have 
\begin{equation}
\sum_{k=1}^{N} \left( - \frac{\vc{q}_k \dpr \nabla T_k}{T_k^2} \right) \geq 0.
\end{equation}
With the above inequalities, it is obvious that the entropy inequality \cref{eq:ent_ineq} holds.
\end{proof}

\section{Numerical method}
\label{sec:numer_meth}
We use the fractional step method to solve the proposed model. The solution procedure is split into the following physical stages:
\begin{enumerate}
\item[(a)] The hydrodynamic sub-system,
\item[(b)] The viscous sub-system,
\item[(c)] The temperature relaxation sub-system,
\item[(d)] The heat conduction sub-system.
\end{enumerate}

The proposed numerical method uses different formulations (the conventional formulation (\ref{eq:reduced_five}) or the symmetric formulation (\ref{eq:reduced_five_1})) for different stages based on the numerical convenience. For hyperbolic stage we use the conventional formulation (\ref{eq:reduced_five}), while for the viscous terms, the thermal relaxation, and the heat conduction we use the corresponding part of the symmetric formulation (\ref{eq:reduced_five_1}).

The viscous and the heat conduction steps contribute (non-linear) parabolic equations with respect to the velocity and the temperature, respectively. Due to the commonality in the algorithm for solving these parabolic PDEs, we separately summarize the corresponding numerical methods in \Cref{subsec:numer_met_para}.

\subsection{Hydrodynamic part}
The hydrodynamic process is governed by the following homogeneous hyperbolic PDEs:  
\begin{subequations} \label{eq:reduced_five_hyper}
\begin{align}
\dudx{\alpha_k \rho_k}{t} + \nabla \dpr \left( {\alpha_k \rho_k \vc{u}} \right) = 0,  \label{eq:reduced_five_hyper:mass} \\
\dudx{\rho\vc{u}}{t} + \nabla\dpr\left(\rho\vc{u}\tpr\vc{u} + p \overline{\overline{I}} \right)  = 0, \label{eq:reduced_five_hyper:mom} \\
 \dudx{\rho E}{t} +  \nabla\dpr \left( \rho E \vc{u} +  p \vc{u} \right) = 0, \label{eq:reduced_five_hyper:en} \\
\dudx{\alpha_l}{t} + \vc{u} \dpr \nabla \alpha_l  =  {\mathcal{R}}_{l}, \label{eq:reduced_five_hyper:vol}
\end{align}
\end{subequations}
where ${\mathcal{R}}_{l} = \alpha_{l} \frac{A-A_{l}}{A_{l}} \nabla \cdot \boldsymbol{u}$.

The hyperbolicity, Riemann invariants and jump conditions of this system in the case of $N=2$ have been investigated in \cite{murrone2005five,kapila2001two} and extensions to the case $N\geq 3$ are straightforward.
Note that the last equation for the volume fraction is in a non-conservative form and is reformulated as 
\begin{equation}\label{eq:reduced_five_hyper:vol1}
\frac{\partial \alpha_{l}}{\partial t}+   \nabla \cdot \left( \alpha_{l} \boldsymbol{u} \right) =  \alpha_l \nabla\cdot \mathbf{u} +  {\mathcal{R}}_{l}.
\end{equation}

The system of \cref{eq:reduced_five_hyper:mass,eq:reduced_five_hyper:mom,eq:reduced_five_hyper:en,eq:reduced_five_hyper:vol1},  can be written in the vector form as follows:
\begin{equation}\label{eq:reduced_five_hyper_conv}
\dudx{\vc{U}}{t} + \nabla\dpr \vc{F}\left(  \vc{U}\right) =  \vc{S}\left(  \vc{U} \right) \nabla\dpr \vc{u},
\end{equation}
where 
\[
\vc{U} = \left[ \alpha_k \rho_k \;\;  \rho u \;\; \rho v \;\; \rho E \;\; \alpha_l  \right]^{\text{T}}, \quad
\vc{F}\left(\vc{U}\right) = u \vc{U} + p \vc{D},
\]
\[
\vc{D}\left(\vc{U}\right) = \left[ 0 \;\;  1 \;\; 0 \;\; u \;\; 0 \right]^{\text{T}},
\quad
\vc{S}\left(\vc{U}\right) = \left[ 0  \;\; 0 \;\;  0  \;\; 0 \;\; \frac{A}{A_l}\alpha_l  \right]^{\text{T}}.
\]

The equation is discretized on a Cartesian grid. 
The conservative part of \cref{eq:reduced_five_hyper_conv} (without the right hand side term) is solved with the Godunov method \cite{Godunov1959}. The numerical flux is determined by using the Riemann solution that is approximated with the HLLC scheme \cite{Toro2009Riemann}. For the high order extension, we adopt the second-order MUSCL scheme and the fifth order WENO scheme \cite{Coralic2014Finite,Johnsen2006Implementation,JIANG1996202} for spatial reconstruction of the local characteristic variables  on cell faces. 

The right hand side term $\vc{S}\left(  \vc{U} \right) \nabla\dpr \vc{u}$ is approximated as
\begin{equation}
\frac{1}{V_{i j k}} \int_{V_{i j k}} \frac{A}{A_l}\alpha_l \nabla \cdot \mathbf{u} \mathrm{d} V \approx \frac{1}{V_{i j k}} \left( \frac{A}{A_l}\alpha_l \right)_{i j k} \int_{{\sigma}_{i j k}} \mathbf{u} \cdot \mathbf{n} \mathrm{d} {\sigma},
\end{equation}
where the subscript $_{ijk}$ denotes the index of the considered cell. ${V_{i j k}}$ and ${{\sigma}_{i j k}}$ are the cell volume and surface, respectively. $\mathbf{n}$ is the surface normal.  The variables $A, \; A_l, \; \alpha_l$ are approximated as the cell-averaged values.

The third-order SSP (Strong Stability-Preserving) Runge–Kutta scheme \cite{Gottlieb2001} is used for time integration. 

The solution of the hyperbolic part is denoted as $\vc{U}_{pri}^{(1)} = [\rho_k^{(1)} \;\;  \vc{u}^{(1)} \;\; p^{(1)} \;\; \alpha_l^{((1))}], $ and serves as the initial data for the following computation.

\subsection{Viscous part}
The governing equations for the viscous step read:
\begin{subequations} \label{eq:reduced_five_viscous}
\begin{align}
\dudx{m_k}{t} = 0,  \label{eq:reduced_five_viscous:mass} \\
\dudx{\rho\vc{u}}{t}  = \nabla \dpr {{\overline{\overline{\tau}}}}, \label{eq:reduced_five_viscous:mom} \\
 \dudx{m_k E_k}{t} = \nabla \dpr ({ \alpha_k \overline{\overline{\tau}}_k \dpr \mathbf{u}}). \label{eq:reduced_five_viscous:enk}
\end{align}
\end{subequations}

The initial data used is the above obtained state $\vc{U}_{pri}^{(1)}$ after the hydrodynamic stage.

In 1D \cref{eq:reduced_five_viscous:mom,eq:reduced_five_viscous:enk} is reduced to the following form:
\begin{subequations} \label{eq:reduced_five_viscous1d}
\begin{align}
 &\rho \frac{\partial u}{\partial t}= \frac{\partial \tau}{\partial x}, \quad \tau = \frac{4}{3} \mu \frac{\partial u}{\partial x}, \label{eq:viscous_para_pde} \\
 &\dudx{m_k E_k}{t} = \frac{\partial \alpha_k \tau_k u}{\partial x}, \quad \tau_k = \frac{4}{3} \mu_k \frac{\partial u}{\partial x}. \label{eq:reduced_five_viscous1d:en}
\end{align}
\end{subequations}

For some application scenarios, the phase viscosity depends on the phase density $\rho_k$ and the temperature $T_k$, i.e., $\mu_k = \mu_k (\rho_k, T_k)$ and   the mixture viscosity $\mu = \mu (\alpha_k, \rho_k, T_k)$. Since the temperature is subject to the impact of the viscosity terms, the parabolic PDE set is non-linear. The numerical methods for solving such parabolic equations are  summarized in \Cref{subsec:numer_met_para}.

Observing \cref{eq:reduced_five_viscous:mass}, it can be seen that the partial densities $m_k$ and the mixture density $\rho$ do not vary at this stage. After solving \cref{eq:reduced_five_viscous}, the variables $E_k$ and $\vc{u}$ are updated. The phase internal energy can be determined as $e_k = E_k - \vc{u} \dpr \vc{u} /2$. 

By using the pressure equilibrium condition
\begin{equation}
p_k \left( \rho_k, e_k \right) = p_k\left( \frac{m_k}{\alpha_k}, e_k \right) = p,
\end{equation}
and the saturation condition for volume fractions $\sum_{k=1}^{N} \alpha_k =1$, we can solve for $\alpha_k$. Having $\alpha_k$ , we can further compute $\rho_k$ and $p$.

The state variable at the end of this stage is denoted as
$
\vc{U}_{pri}^{(2)}.
$

\subsection{Temperature relaxation part}\label{subsec:TR}
The governing equations for the temperature relaxation are as follows:
\begin{subequations} \label{eq:reduced_five_tr}
\begin{align}
\dudx{m_k}{t} &= 0,  \label{eq:reduced_five_tr:mass_1} \\
\dudx{\rho\vc{u}}{t}  &= 0, \label{eq:reduced_five_tr:mom_1} \\
\dudx{m_k E_k}{t} &= \mathcal{Q}_k^{sum,1}, \label{eq:reduced_five_tr:enk_1} 
\end{align}
\end{subequations}
with $\vc{U}_{pri}^{(2)}$ as the initial data. $\mathcal{Q}_k^{sum,1}$ is the heat transfer at this stage, which drives the non-equilibrium temperatures into an equilibrium state.

Simple algebraic manipulations of  \cref{eq:reduced_five_tr} give
\begin{equation}\label{eq:mkdek}
m_k^{(3)} = m_k^{(2)}, \;\; {\vc{u}}^{(3)} = {\vc{u}}^{(2)}, \;\; m_k \dudx{e_k}{t} = \mathcal{Q}_k^{sum,1},
\end{equation}
where the superscript ${(3)}$ represents the variables after the temperature relaxation. It can be seen that the partial density $m_k$ and the velocity $\vc{u}$ do not vary at this stage.

We solve the governing equations for the temperature relaxation maintaining pressure equilibrium:
 \begin{equation}\label{eq:pres_eq}
 \dudx{p_k}{t} = \dudx{p}{t} \quad \text{or} \quad p_k\left( T_k, \rho_k \right) = p.
 \end{equation}

The saturation condition for volume fractions leads to
\begin{equation}\label{eq:vol_saturation}
\sum_{k=1}^{N} \frac{m_k}{\rho_k} = 1.
\end{equation}

With \cref{eq:pres_eq,eq:vol_saturation}, the phase density can be expressed as 
\begin{equation}\label{eq:rhok_fun}
\rho_k = \rho_k \left(m_1, m_2, \cdots, m_N, T_1, T_2, \cdots, T_N \right).
\end{equation}

Moreover,
\begin{equation}\label{eq:ek_fun}
e_k = e_k \left( \rho_k, T_k \right) = e_k \left(m_1, m_2, \cdots, m_N, T_1, T_2, \cdots, T_N \right).
\end{equation}

Combination of \cref{eq:mkdek,eq:reduced_five_tr,eq:ek_fun} leads to
\begin{align}\label{eq:tr_eq1}
\sum_{j=1}^{N}\mathcal{A}_{jk}\dudx{T_j}{t}  = \mathcal{Q}_k^{sum,1}, 
\end{align}
where \[\mathcal{A}_{jk} = m_k \dudx{e_k}{T_j}.\]

The time derivative is approximated as
\begin{equation}\label{eq:dTk_dis}
\dudx{T_k}{t} = \frac{T_k^{(3)} - T_k^{(2)}}{\Delta t}.
\end{equation}

An instantaneous temperature relaxation is assumed within the hydrodynamic time step $\Delta t$, thus, we have 
\begin{equation}\label{eq:Teq}
T_1^{(3)} =T_2^{(3)} =\cdots =T_N^{(3)} = T^{(3)}.
\end{equation}

Combining \cref{eq:tr_eq1,eq:dTk_dis,eq:Teq} and having in mind $\sum_{k=1}^{N} \mathcal{Q}_{k}^{sum,1} = 0$, one can obtain:
 \begin{equation}\label{eq:tr_temp_av}
 T^{(3)} = \frac{ \sum_{k=1}^{N} \sum_{j=1}^{N} \mathcal{A}_{jk} T_j^{(2)}}{\mathcal{A}},
 \end{equation}
 where $\mathcal{A} = \sum_{k=1}^{N} \sum_{j=1}^{N} \mathcal{A}_{jk}.$
 
Having $T^{(3)}$, we can solve for $\rho_k^{(3)}$ with \cref{eq:rhok_fun}, and then for $p^{(3)}$ with \cref{eq:pres_eq}. Since the partial density does not vary, i.e. $m_k^{(3)} = m_k^{(2)}$, the volume fractions can be evaluated with $\alpha_l^{(3)} = m_l^{(3)}/\rho_l^{(3)}$. In this way, we can determine the temperature-relaxed state in each cell as
\begin{equation}\label{eq:Upri}
\vc{U}_{pri}^{(3)} = [\rho_k^{(3)} \;\; \vc{u}^{(3)} \;\; p^{(3)} \;\; \alpha_l^{((3))}] .
\end{equation}

\subsection{Heat conduction part}\label{subsec:HC}
The governing equations for the heat conduction read:
\begin{subequations} \label{eq:reduced_five_hc}
\begin{align}
\dudx{m_k}{t} &= 0,  \label{eq:reduced_five_hc:mass_1} \\
\dudx{\rho\vc{u}}{t}  &= 0, \label{eq:reduced_five_hc:mom_1} \\
\dudx{m_k E_k}{t} &=  q_k + \mathcal{I}_k + \mathcal{Q}_k^{sum,2}, \label{eq:reduced_five_hc:enk_1} 
\end{align}
\end{subequations}
where $\mathcal{Q}_k^{sum,2}$ is the inter-phase heat conduction that holds the temperature in equilibrium at this stage.

The initial condition used is the state after the temperature relaxation, $\vc{U}_{pri}^{(3)}$ in \cref{eq:Upri}. From \cref{eq:reduced_five_hc}, one can deduce
\begin{equation}
m_k^{(4)} = m_k^{(3)}, \;\; {\vc{u}}^{(4)} = {\vc{u}}^{(3)}, \;\; m_k \dudx{e_k}{t} = q_k + \mathcal{I}_k + \mathcal{Q}_k^{sum,2}.
\end{equation}

We solve the governing equations for  the heat conduction maintaining the temperature and pressure equilibrium:
  \begin{equation}\label{eq:PT_condition}
 \dudx{T_k}{t} = \dudx{T}{t}, \quad \dudx{p_k}{t} = \dudx{p}{t}. 
 \end{equation}

Since the pressure equilibrium is maintained, \cref{eq:rhok_fun} still holds. With \cref{eq:reduced_five_hc,eq:rhok_fun}, one can obtain:
\begin{align}\label{eq:hc_eq1}
\sum_{j=1}^{N} \mathcal{A}_{j k} \frac{\partial T_{j}}{\partial t} = q_k + \mathcal{I}_k + \mathcal{Q}_k^{sum,2},
\end{align}

At this stage the phase temperatures are initially relaxed, i.e., $T_k^{(3)} = T^{(3)}$. Moreover, their time derivatives are maintained to be in equilibrium due to \cref{eq:PT_condition}. This derivative equilibrium condition can be achieved by properly defining the inter-phase heat conduction term $\mathcal{Q}_k^{sum,2}$. As a result, the phase temperatures are in equilibrium through the heat conduction stage, i.e., $T_k(t) = T(t)$ at a fixed spatial location.

Summing \cref{eq:hc_eq1} over $k$ and having in mind $\sum_{k=1}^{N} \mathcal{Q}_k^{sum,2} = 0$ yield
\begin{align}\label{eq:hc_para_pde}
\mathcal{A} \dudx{T}{t} = \sum_k^{N} \left( q_k + \mathcal{I}_k \right) = q + I,
\end{align}
where  $q = \nabla \dpr \left( \lambda \nabla T \right)$, and the mixture heat conduction coefficient is defined as $\lambda = \sum_{k=1}^{N} \alpha_k \lambda_k$, $I = \sum_{k=1}^{N} \mathcal{I}_k$.

Solving the nonlinear parabolic equation (\ref{eq:hc_para_pde}) in the computational domain with the initial data $T^{(3)}$, we can determine the temperature after the heat conduction stage -- $T^{(4)}$. Having $T^{(4)}$, we then determine the corresponding state variables in the way similar to that of the temperature relaxation in \cref{subsec:TR}:  
\begin{equation}\label{eq:Upri2}
\vc{U}_{pri}^{(4)} = [\rho_k^{(4)}  \;\; \vc{u}^{(4)} \;\; p^{(4)} \;\; \alpha_l^{((4))}] .
\end{equation}

\subsection{Numerical methods for the parabolic diffusion PDEs}\label{subsec:numer_met_para}
The parabolic PDEs (\cref{eq:viscous_para_pde,eq:hc_para_pde}) in the above steps have the following general form (in 1D):
\begin{equation}\label{eq:nonlin_para_pde}
\frac{\partial v}{\partial t}=L [v]+f(x, t), \quad x \in \Lambda \subset \mathbb{R}, 
\end{equation}
where $L[\cdot]$ is a quasi-linear elliptic positive definite operator.
\begin{equation}
L[v] = \dudx{}{x}\left( k(v) \dudx{v}{x} \right).
\end{equation}


We adopt the simple iterative method to solve this non-linear PDE, which produces the following iteration sequences
\begin{equation}\label{eq:lin_para_pde}
\frac{\partial v^{(s+1)}}{\partial t} = \dudx{}{x}\left( k(v^{(s)}) \dudx{v^{(s+1)}}{x} \right) +f(x, t),
\end{equation}
where the superscript $s$ denotes the iteration number and $v^{(0)}$ is taken to be the values at the beginning of the current computation stage, i.e., $v^{(0)} = v^n$. The non-linear coefficient $k$ may also depend on the other state variables that are taken to be those from the last iteration.

In each iteration we solve the linearized parabolic \cref{eq:lin_para_pde} with the traditional explicit or implicit schemes with respect to $v^{(s+1)}$. The iteration is performed until $||v^{(s+1)} - v^{(s)}||$ is less than a prescribed threshold. In fact, the discretization form of \cref{eq:lin_para_pde} can also be viewed as a discretization of \cref{eq:nonlin_para_pde} with specific linearization.

In the present work we implement both the implicit and explicit schemes. In the former the linear equations are solved with the preconditioned conjugate gradient (PCG) method. 
The latter is realized with the monotonicity-preserving  local iteration method (LIM) \cite{Zhukov2010} based on the Chebyshev parameters. Each iteration of the LIM is equivalent to the implementation of a traditional 7-point (in 3D) explicit scheme, but allows a much larger stable time step, which is of the order $\mathcal{O}(P^2 \Delta x^2)$ ($P$ is the number of stencil points, $\Delta x$ is the mesh size). Its parallel implementation is as straightforward as the traditional explicit scheme and demonstrates good scalability and efficiency \cite{zhukov2018}.

For simplicity, we demonstrate the LIM solution procedure in 1D on a  Cartesian grid with uniform spacing $\Delta x$. In each iteration, we need to solve the following linear parabolic equation
\begin{equation}\label{eq:lin_para_pde1}
\frac{\partial \overline{v}}{\partial t} = {L^{\prime}}[\overline{v}] +f(x, t), \quad  {L^{\prime}}[\overline{v}] = \dudx{}{x}\left( \uline{k}  \dudx{\overline{v}}{x} \right)
\end{equation}
where $\overline{v} = v^{(s+1)}$, $\uline{k} = k(v^{(s)})$. The operator ${L^{\prime}}$ is approximated with a second-order central difference operator ${L_h^{\prime}}$, which is self-adjoint and has real positive eigenvalues within an interval $[\lambda_{\text{min}}, \lambda_{\text{max}}]$,
\begin{equation}
 L_h^{\prime}[{{\overline{v}}}] = 
 \frac{{\uline{k}}_{i+1/2} \left( {\overline{v}}_{i+1} - {\overline{v}}_{i} \right) - {\uline{k}}_{i-1/2} \left( {\overline{v}}_{i} - {\overline{v}}_{i-1} \right)}{\Delta x^2}.
\end{equation}

The LIM can be written in the form of the following $2P-1$ iterations
\begin{equation}\label{eq:wmm}
{\overline{v}}^{(m)}=\frac{1}{1+\Delta t b_{m}}\left({\overline{v}}^{n}+\Delta t b_{m} {\overline{v}}^{(m-1)} + \Delta t L_{h} {\overline{v}}^{(m-1)}+\Delta t f^{(n)}\right),
\end{equation}
where $m=1,2, \cdots, 2 P-1$. The stencil size depends on the time step $\Delta t$,  $P= \lceil\pi / 4 \sqrt{\Delta t \lambda_{\max }+1}\rceil$, $\lceil x \rceil = \text{min}\{n\in \mathbb{Z} \; | \; n\geq x\}$. $b_m$ is a set of iteration parameters
\begin{equation}
\left(b_{1}, b_{2}, \cdots, b_{2 P-1}\right)=\left(a_{P}, a_{P-1}, \cdots, a_{2}, a_{P}, a_{P-1}, \cdots, a_{1}\right), \quad a_{m}=\frac{\lambda_{\max }}{1+\beta_{1}}\left(\beta_{1}-\beta_{m}\right),
\end{equation}
where $\beta_{m}$ is the root of the  Chebyshev polynomial $T_P(x)$: $\text{cos}[{(2m-1)\pi}/{(2P)}]$.

The last iteration of \cref{eq:wmm} can be written as
\begin{equation}
{\overline{v}}^{(2P-1)}= \left({\overline{v}}^{n} + \Delta t L_{h} {\overline{v}}^{(2P-2)} + \Delta t f^{(n)}\right),
\end{equation}
here, we see that ${\overline{v}}^{(2P-2)}$ can be viewed as a predicted solution.

For details of LIM, see \cite{Zhukov2010,zhukov2018}. 

\subsubsection{LIM for the viscous part}
When LIM is used for solving the viscous part of the parabolic PDE \cref{eq:viscous_para_pde}, the total energy and volume fraction should be updated accordingly.  

The last LIM iteration for solving \cref{eq:viscous_para_pde} is 
\begin{equation}\label{eq:dis_vis_mom}
\rho_i \frac{u^{n+1}_i - u^{n}_i}{\Delta t} = \frac{\tau_{i+1/2} - \tau_{i-1/2}}{\Delta x},
\end{equation} 
\[
\tau_{i+1/2} = \frac{4}{3} \mu_{i+1/2} \left( \dudx{u}{x} \right)_{i+1/2}, \quad \left( \dudx{u}{x} \right)_{i+1/2} = \frac{u_{i+1}^{(2P-2)} - u_i^{(2P-2)}}{\Delta x},
\]
\[\mu_{i+1/2} = \sum_{k=1}^{N} (\alpha_k)_{i+1/2} \mu_k, \quad (\alpha_k)_{i+1/2} = [(\alpha_k)_{i} +(\alpha_k)_{i+1}]/2,\]
the subscripts $i+1/2$ and $i-1/2$ denote the values on the left and right cell faces, respectively.

The total energy is updated as 
\begin{equation}\label{eq:dis_vis_en}
m_k \frac{\left(E_k\right)^{n+1}_{i}-\left(E_k\right)^{n}_{i}}{\Delta t} = \frac{ (\alpha_k \tau_k  u)_{i+1/2} -  (\alpha_k \tau_k u)_{i-1/2}}{\Delta x}
\end{equation}
where  
\[ \quad u_{i+1/2} = (u_{i+1}^{(2P-2)} + u_i^{(2P-2)})/2, \quad \left(\tau_k\right)_{i+1/2} = \frac{4}{3} (\mu_k)_{i+1/2} \left( \dudx{u}{x} \right)_{i+1/2}.\]
Note that \cref{eq:dis_vis_en} is solved only in the last LIM iteration  and the velocities used here are the predicted ones after $2P-2$ LIM iterations, i.e., $ u^{(2P-2)}$.

The above algorithm can be extended to multiple dimensions with the dimension-splitting technique and to non-uniform grids with the corresponding schemes to approximate the variables and gradients on cell faces.

\subsubsection{LIM for the heat conduction part}
The algorithm for solving the heat conduction equation (\ref{eq:hc_para_pde}) is similar to that for the viscous part. Note that here a source term exists, i.e.,
 \[ f = \sum_{k=1}^{N} \mathcal{I}_k.\]
This term is approximated in an explicit way by directly setting it to be the computed values from the last time step.

\subsection{Preservation of the PVT equilibrium condition}

In \cite{Alahyari2015,Johnsen2012Preventing} two different mixture EOSs are proposed to maintain the PVT property, which leads to the ambiguity in the definition of mixture EOS and other thermodynamic variables (e.g., the mixture entropy). Thus, the consistency with the second law of thermodynamics is an issue of much controversy. The proposed model formulation ensures the PVT property and the thermodynamical laws with a uniquely defined EOS. 

Let us consider the following Riemann problem:
\begin{align}\label{eq:initial1}
u^{L} =u^{R}=u > 0, \;\;
\rho_{k}^{L} =\rho_{k}^{R}=\rho_{k},  \;\;
e_{k}^{L} = e_{k}^{R} = e_{k}, \;\;
{\alpha_l}^{L} \neq {\alpha_l}^{R}, \;\;  p_{L} = p_{R} = p, \;\;  T_{L} = T_{R} = T.
\end{align}

Here, the phase temperatures are initially in equilibrium, i.e., 
\begin{equation}\label{eq:initial2}
T_{k,L} = T_{L}, \; T_{k,R} = T_{R}.
\end{equation}

\begin{proposition}
The solution to the proposed model (\ref{eq:reduced_five}) preserves the PVT property with the initial discontinuity  (\ref{eq:initial1}) and (\ref{eq:initial2}).
\end{proposition}

\begin{proof}\label{proof:1}
The model equations are solved in the framework or the Godunov finite volume method. We are concerned about the state evolution of the cell downstream the initial discontinuity $\vc{U}^{*}$ after one time step.

When the velocity is divergence free, the model equations (\ref{eq:reduced_five}) without diffusions or external energy source is reduced to that in \cite{allaire2002five}. Therefore, we borrow the well proved proposition on pressure equilibrium formulated in \cite{allaire2002five,Zhangchao2020}:
\begin{equation}
u^{*} = u, \quad
\rho_{k}^{*} =\rho_{k}, \quad
e_{k}^{*} =e_{k}, \quad
p^{*} =p.
\end{equation}
Since $T_k^{*} = T_k (\rho_k^{*}, p^{*})$, one can obtain 
\begin{equation}
T_k^{*} = T.
\end{equation}
This means that the hydrodynamic step does not undermine the PVT property.

The temperature relaxation procedure \cref{eq:tr_temp_av} can be viewed as a specific averaging of the phase temperatures and thus does not alter the equilibrium temperature. Further, according to \cref{eq:rhok_fun}, phase densities $\rho_k$ and volume fractions $\alpha_k$ remain unchanged during the temperature relaxation since phase temperatures $T_k$ and partial densities $m_k$ do not vary. As a function of the phase temperature and the phase density, the pressure does not vary at this stage either.

Since temperature and velocity are spatially uniform, the thermal conductive and viscous terms have no impact on the solution.
\end{proof}

\section{Numerical results}
\label{sec:numer_res}
In this section, we perform several numerical tests to verify the proposed model and numerical methods. The considered tests demonstrate the impact of temperature relaxation,  viscosity and heat conduction on the numerical solutions. The numerical results with heat conduction of the proposed reduced model is compared to those of the one-temperature model \cite{Alahyari2015}, demonstrating the superiority of the present model in convergence performance, especially when complicated EOS is involved. We also apply the model for simulating the laser ablation of a multicomponent target in the ICF field.

\subsection{The pure translation of a material interface}
In this section we verify the capability of the proposed model to maintain the PVT property. 
The computational domain is $[0.0\text{m},1.0\text{m}]$.
Two fluids are separated by the material interface located at $x = 0.2$m. The properties of the two fluids on the left and right of the interface are as follows:
\[\gamma_1 = 4.40, \; p_{\infty, 1} = 6.00 \times 10^{6} \text{Pa}, \; C_{v,1} = 58.82\text{J}/(\text{kg}\cdot \text{K}), \; \alpha_1 = 1 - \varepsilon, \]
\[\gamma_2 = 1.40, \; p_{\infty, 2} = 0.00 \text{Pa}, \; C_{v,2} = 125.00\text{J}/(\text{kg}\cdot \text{K}), \; \alpha_1 = \varepsilon, \]
where $\varepsilon$ is a small positive number, and is taken as $1.00 \times 10^{-6}$ if not mentioned specifically.

The initial pressure, temperature and velocity are uniformly distributed and are $1.00 \times 10^5$Pa, $3.00\times 10^3$K, $1.00\times 10^2$m/s, respectively. Phase densities are determined via their EOSs. Extrapolation boundary conditions are imposed on both sides of the computational domain. 

Computations last to the time moment $t = 5.00 \times 10^{-6}$s. 
We solve this problem with two models: the fully conservative one-temperature four-equation model \cite{Lemartelot2014} and the proposed reduced model.  
Here, the second order MUSCL scheme with MINMOD limiter is used. 

\begin{figure}[h!]
\centering
\subfloat[Pressure]{\label{figPA:dens}\includegraphics[width=0.5\textwidth]{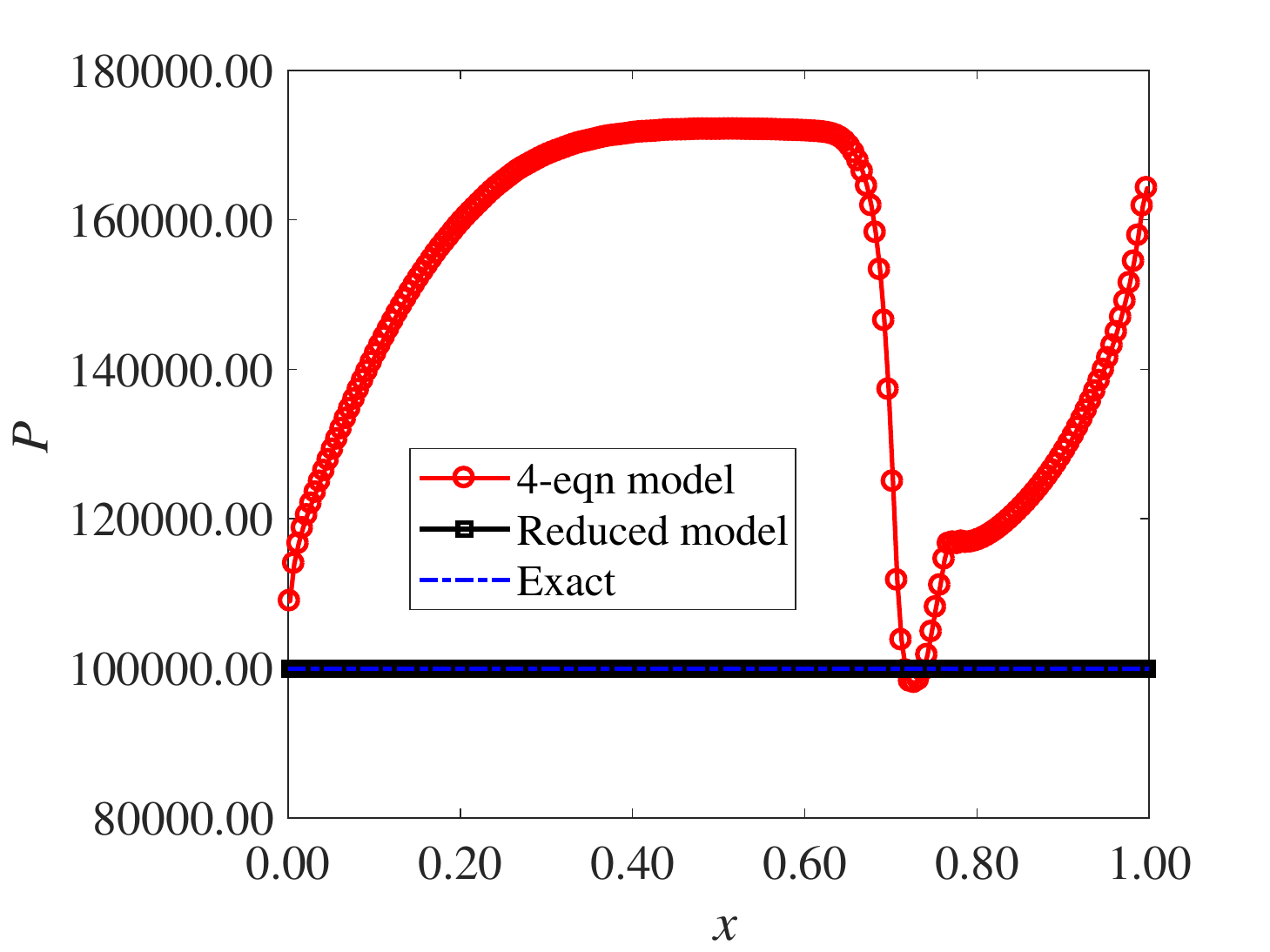}}
\subfloat[Velocity]{\label{figPA:dens1}\includegraphics[width=0.5\textwidth]{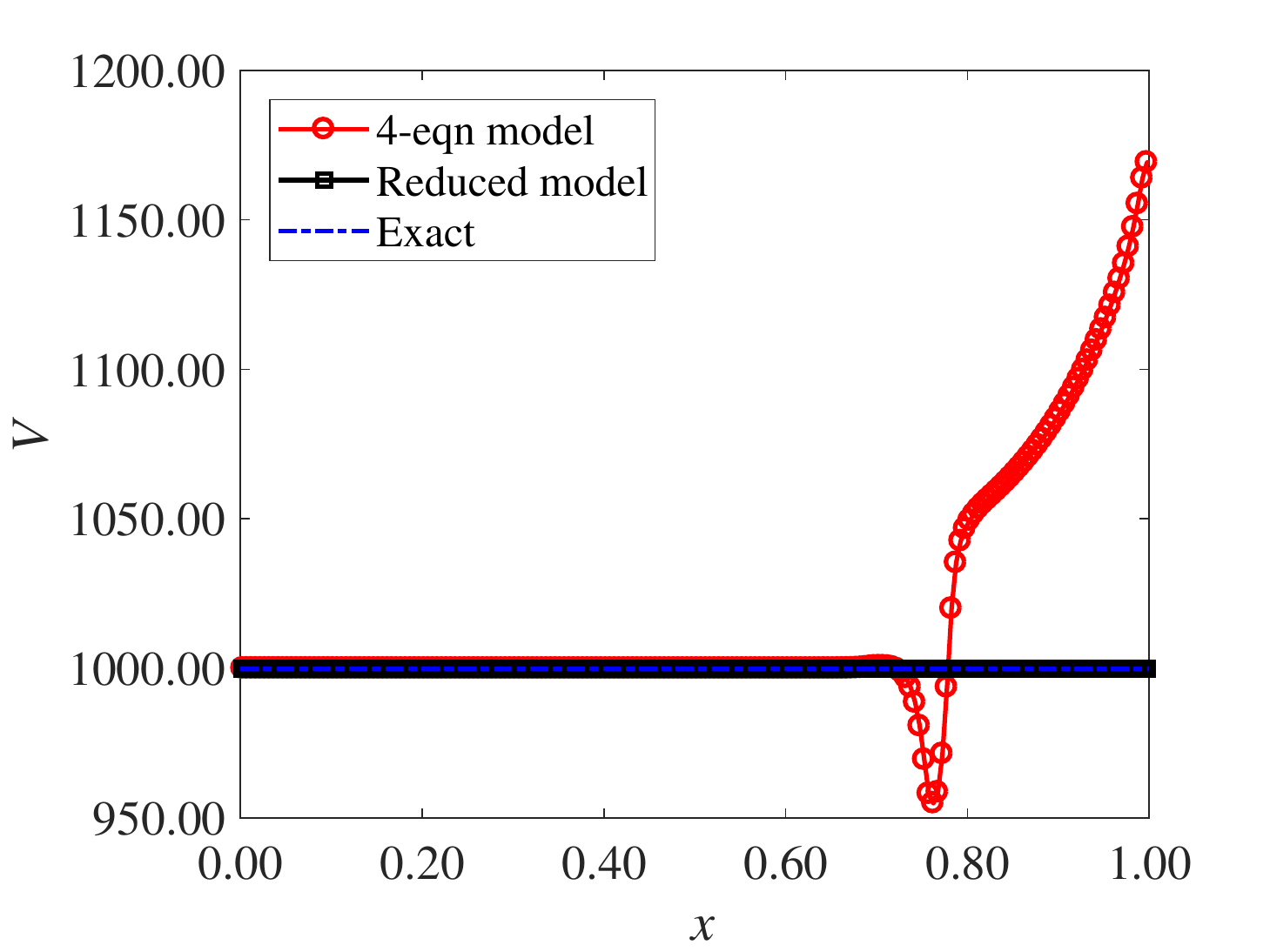}}\\
\subfloat[Temperature]{\label{figPA:tem}\includegraphics[width=0.5\textwidth]{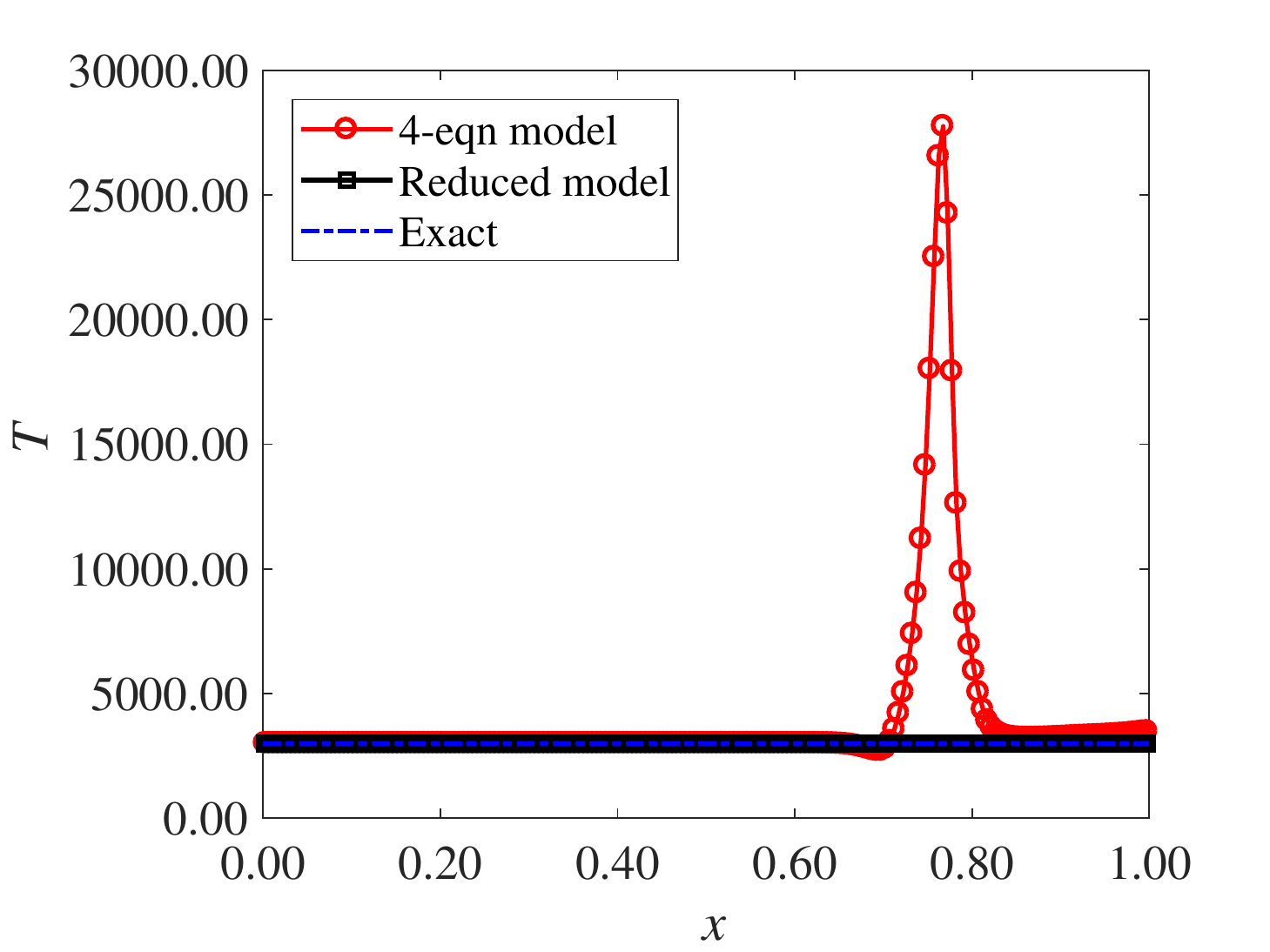}} 
\subfloat[Density]{\label{figPA:tem1}\includegraphics[width=0.5\textwidth]{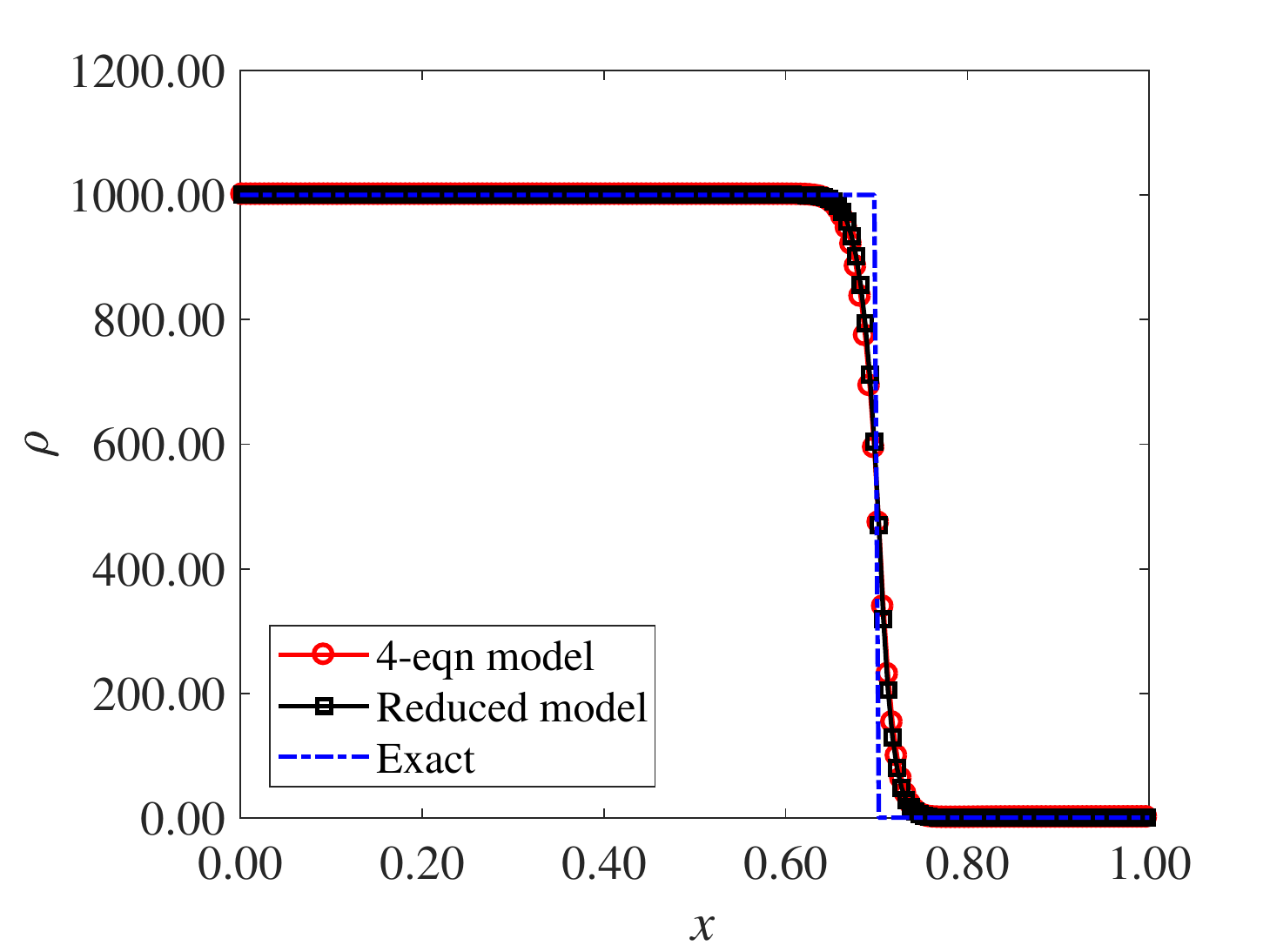}}
\caption{Numerical results for the pure translation problem. Red circles -- results of the conservative four-equation model, black squares -- results of the proposed reduce model.}
\label{fig:figPA} 
\end{figure}

For this problem, the pressure, the velocity and the temperature should remain unchanged in the exact solution. The numerical results of the considered models are compared with the exact ones in \Cref{fig:figPA}. One can see that the four-equation model results in spurious oscillations in pressure, velocity and temperature. The proposed reduced model maintains the PVT property as expected.

\subsection{The multicomponent shock tube problem}
In this section we consider a two-fluid shock tube problem with heat conduction. Two fluids are initially at rest and and separated by the material interface located at $x = 0.7$m separating them. The parameters $\gamma$ and $p_{\infty}$ of the two fluids are the same as the last test and $C_{v,1} = 1606.00\text{J}/(\text{kg}\cdot \text{K}), \;\; C_{v,2} = 714.00\text{J}/(\text{kg}\cdot \text{K})$.
The pressure and the temperature on the left and right of the interface are 
 $p = 1.00\times 10^9, \; T = 293.02\text{K}$
and $p = 1.00\times 10^5, \; T = 7.02\text{K}$, respectively.
Phase densities are determined in such a way that the initial temperatures are in equilibrium.

The heat conduction coefficients are set to be $\lambda_1 =  1.00 \times 10^4 {\text{W}/(\text{m}\cdot\text{K})}$ and $\lambda_2 =  1.00 \times 10^6 {\text{W}/(\text{m}\cdot\text{K})}$, respectively, which are large enough to demonstrate the impact of the heat conduction on the numerical results.

\begin{figure}[h!]
\centering
\subfloat[Temperature]{\label{fig00:temps}\includegraphics[width=0.5\textwidth]{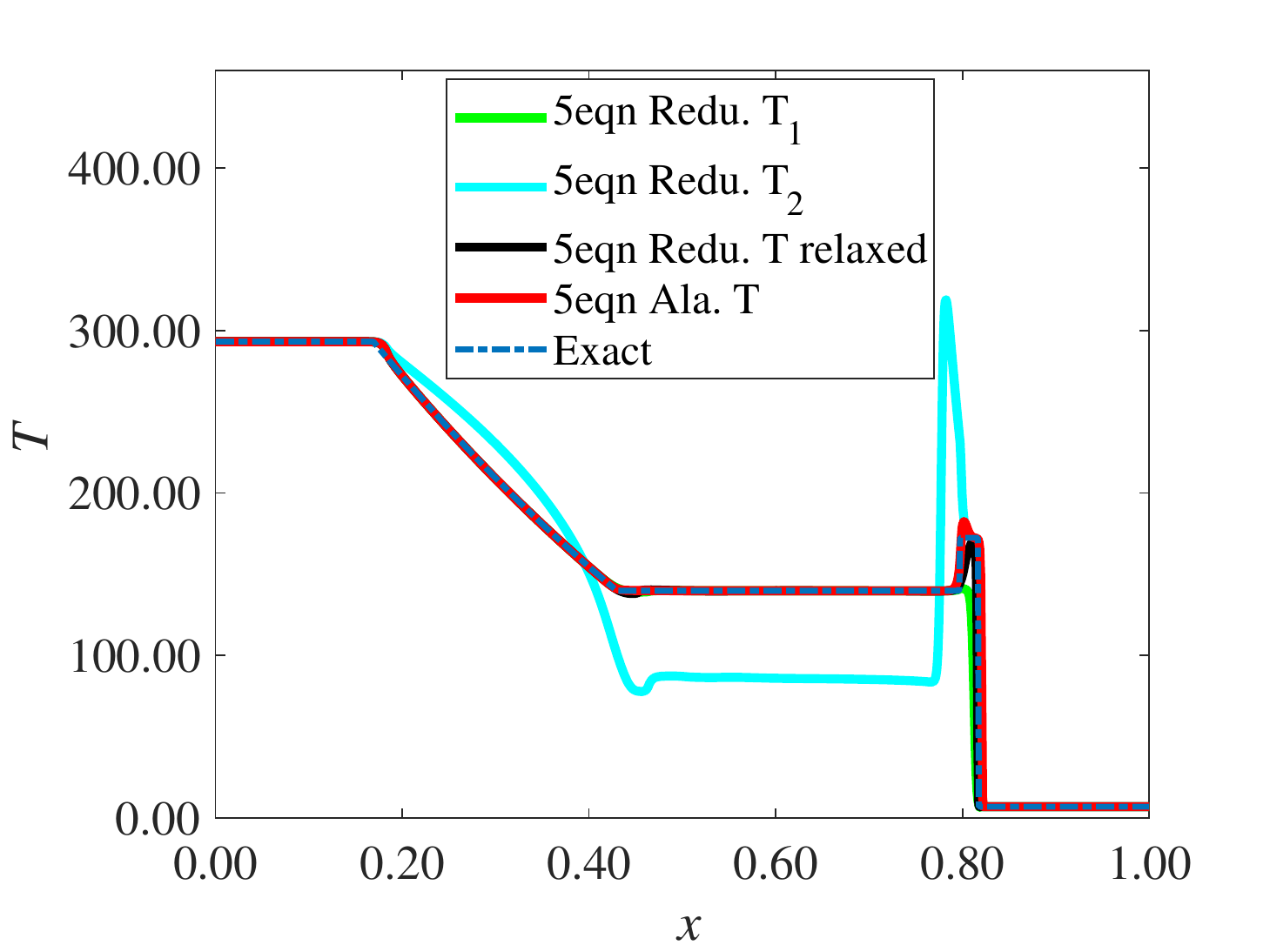}}
\subfloat[Temperature, locally enlarged]{\label{fig00:temps1}\includegraphics[width=0.5\textwidth]{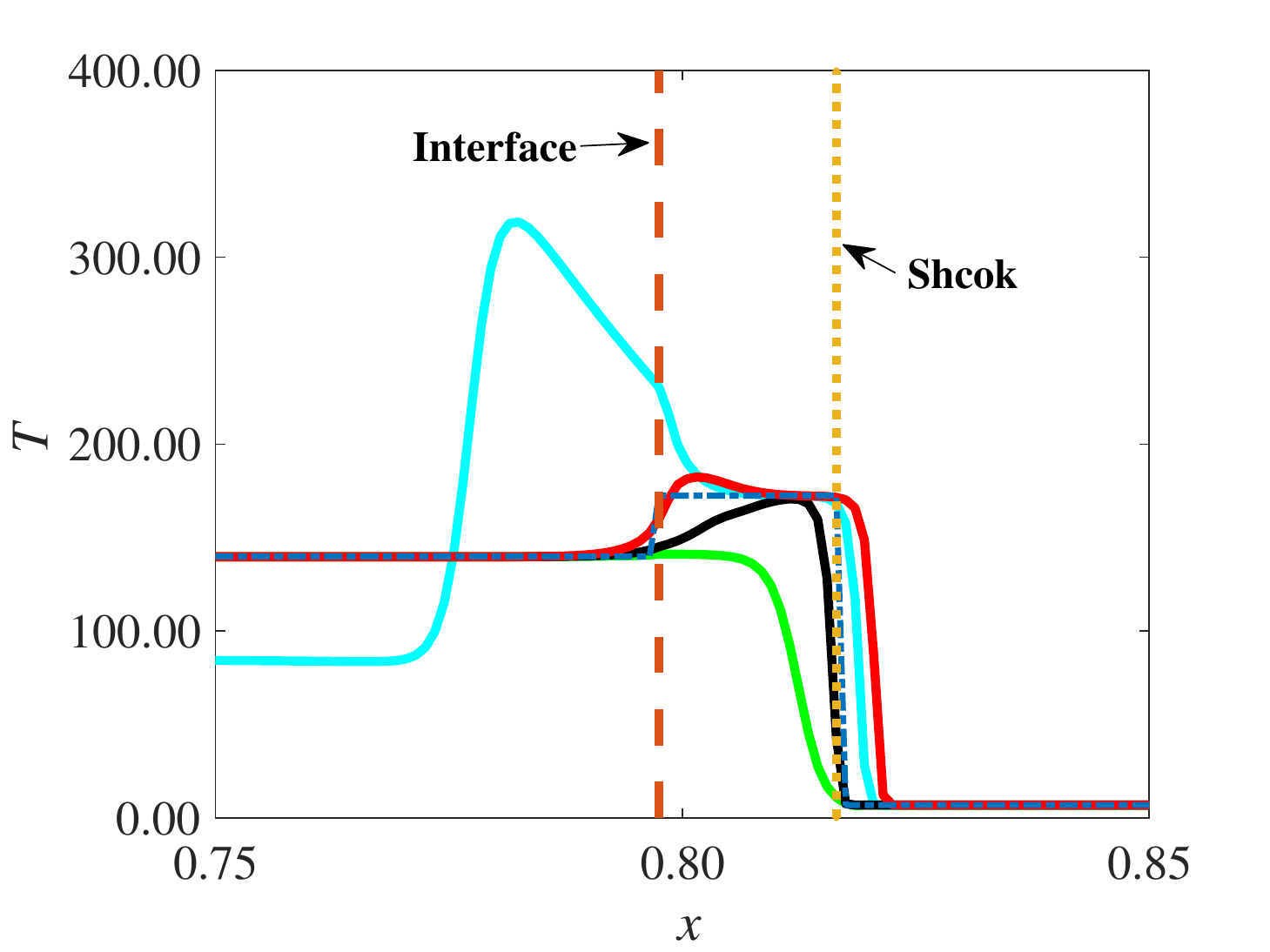}}
\caption{Temperatures obtained with different models after the hydrodynamic step.}
\label{fig:fig00} 
\end{figure}

\begin{figure}[h!]
\centering
\subfloat[Density]{\label{fig:dens}\includegraphics[width=0.5\textwidth]{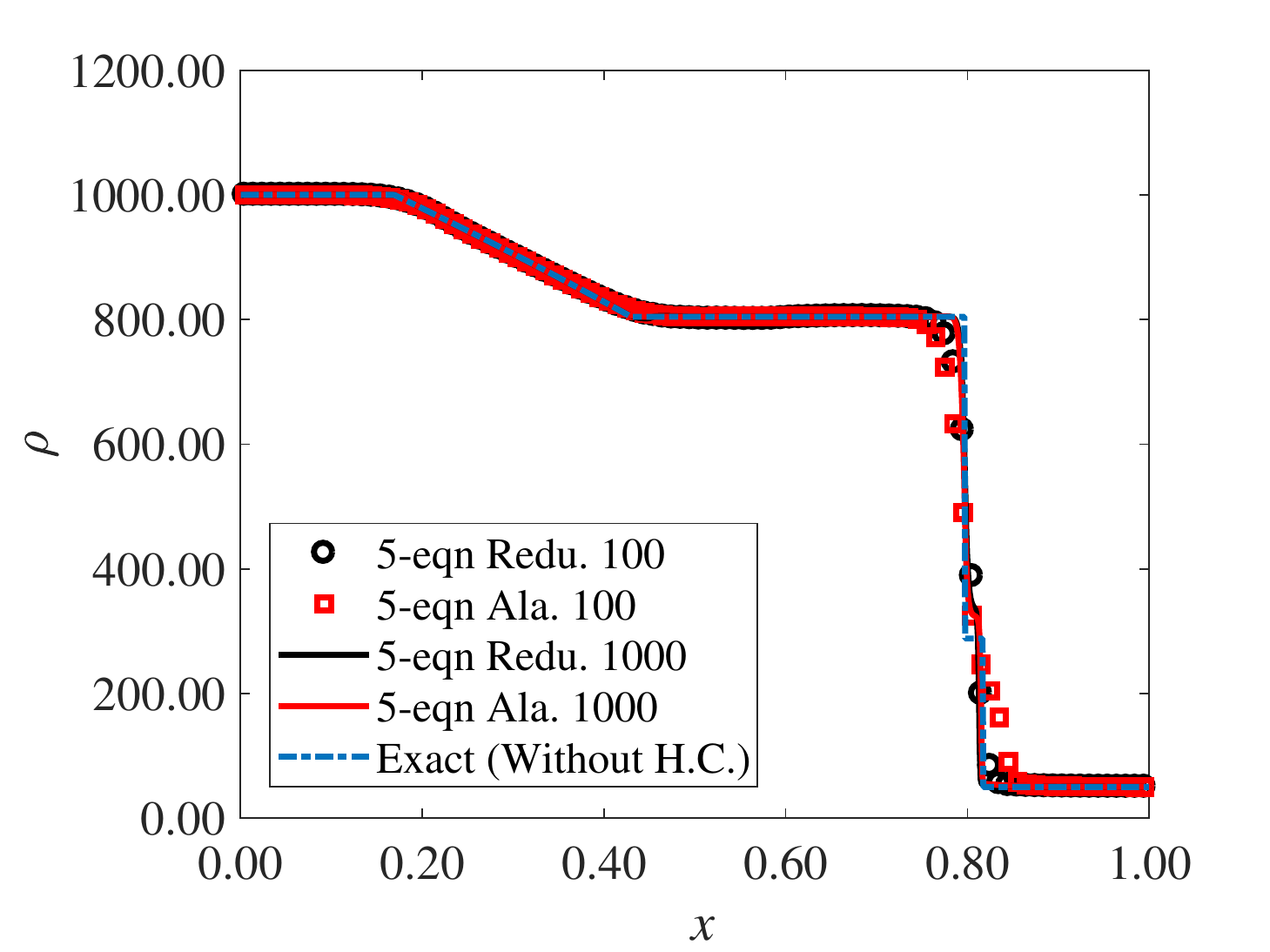}}
\subfloat[Density, locally enlarged]{\label{fig:dens1}\includegraphics[width=0.5\textwidth]{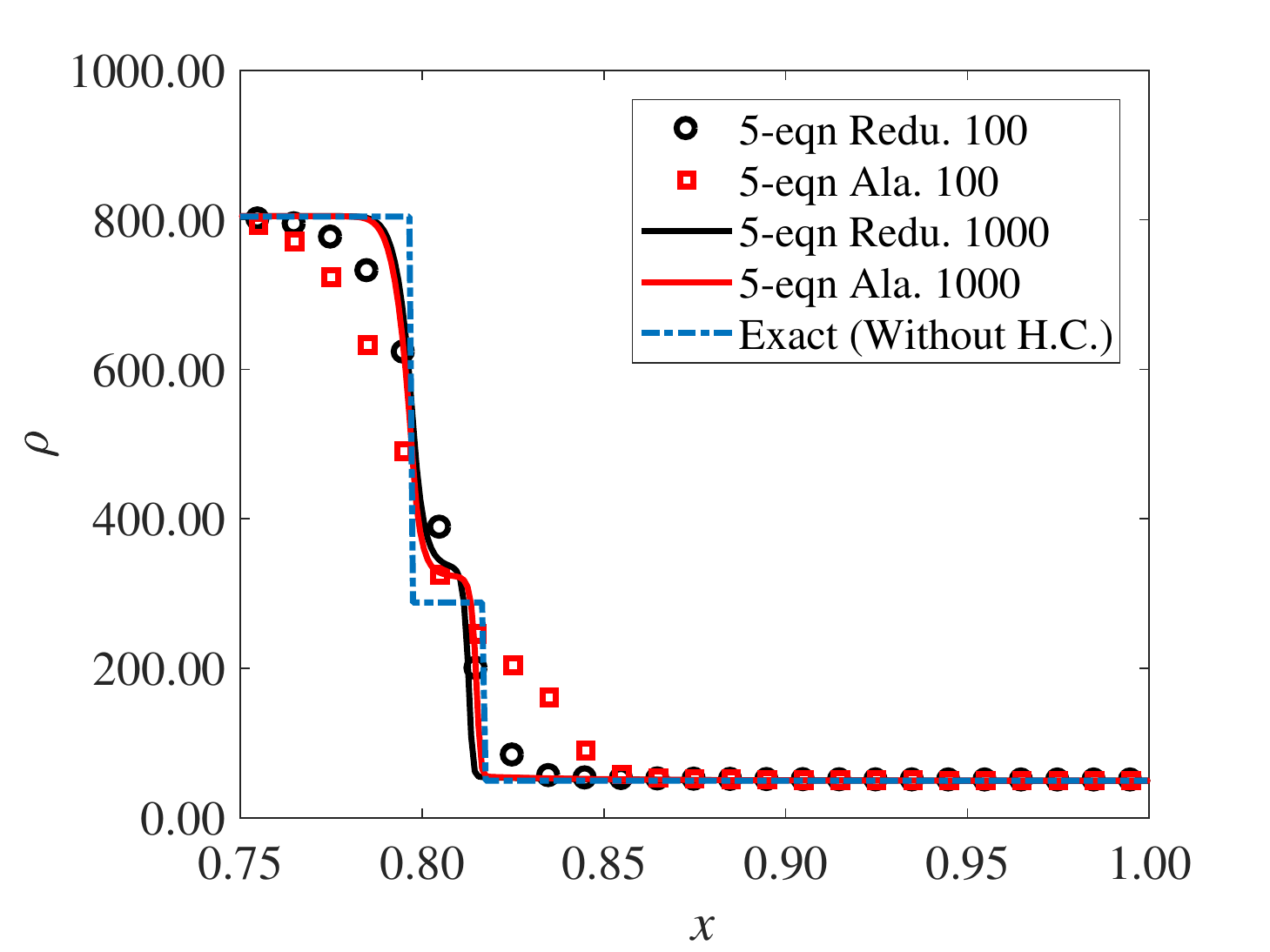}}\\
\subfloat[Temperature]{\label{fig:tem}\includegraphics[width=0.5\textwidth]{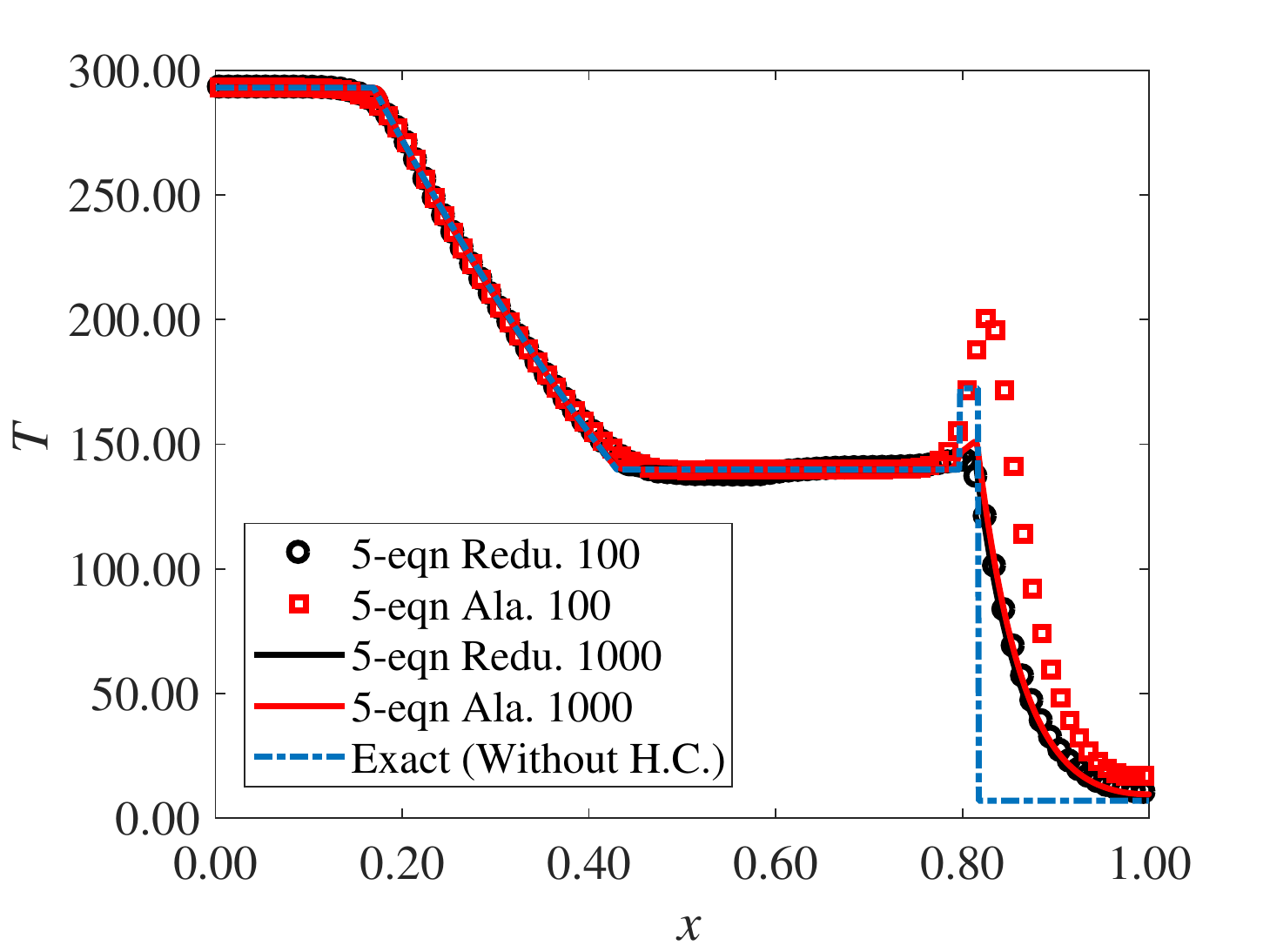}} 
\subfloat[Temperature, locally enlarged]{\label{fig:tem1}\includegraphics[width=0.5\textwidth]{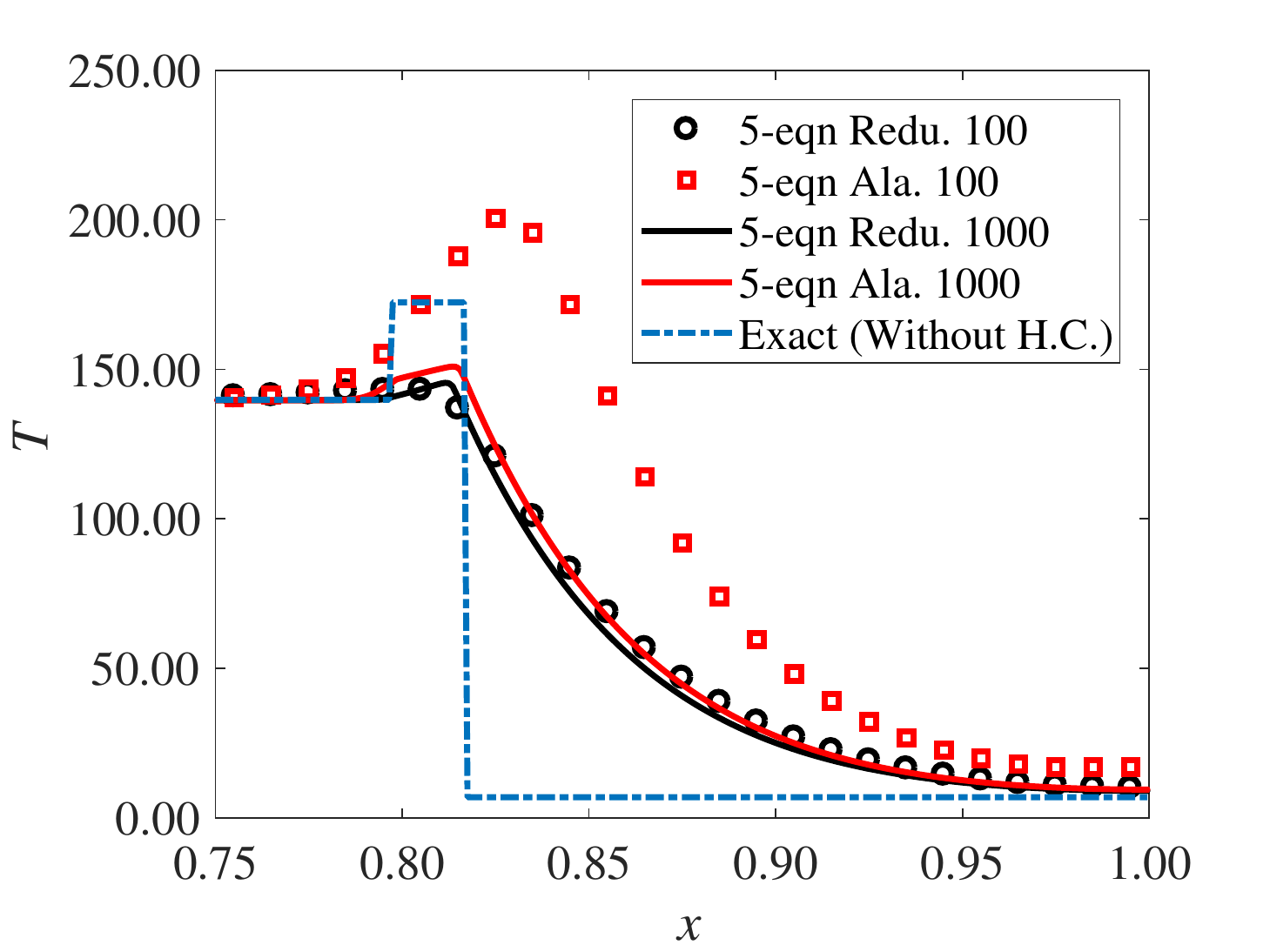}}
\caption{Numerical results for the density, the temperature and the velocity in the water-gas shock tube problem.}
\label{fig:fig01} 
\end{figure}

The numerical results obtained on the 100-cell and 1000-cell grids with different models. Here the second-order MUSCL scheme is used.
We first perform computations without any diffusion processes.  The numerical results for temperature are displayed in \Cref{fig:fig00}. The phases temperatures ($T_1$ and $T_2$ in the figure) are obtained by solving the reduced model without the temperature relaxation. The relaxed temperature is obtained by solving the reduced model with the temperature relaxation. We compare the numerical results of the reduced model to those of the one-temperature five-equation model of Alahyari Beig \cite{Alahyari2015,Johnsen2012Preventing} that are denoted as ``5-eqn Ala. T'' in the figure. The exact solution displayed is the analytical solution to the Riemann problem without any relaxations or diffusions. The phase temperatures $T_1$ and $T_2$ are only meaningful in the domain occupied by the corresponding phase, i.e., the phase temperature $T_1$ ($T_2$) is only meaningful to the left (right) of the material interface, as demonstrated in \Cref{fig00:temps1}. The temperature relaxation drives these phase temperatures in an equilibrium temperature (the line denoted as ``5eqn Redu. T relaxed'' in the figure). 
Moreover, the one-temperature model tends to overestimate the shock speed in the cases without the heat conduction (\Cref{fig00:temps,fig00:temps1}) or with (\Cref{fig:dens,fig:dens1}).

From \Cref{fig:tem,fig:tem1} it can be seen that the solution of both models tend to converge to the same solution. 
However, the convergence performance of the proposed temperature non-equilibrium model is evidently superior to that of the one-temperature one. We attribute the advantage in numerical behaviour of the proposed model to its consistency  with the second law of thermodynamics.

\begin{figure}[h!]
\centering
{\includegraphics[width=0.5\textwidth]{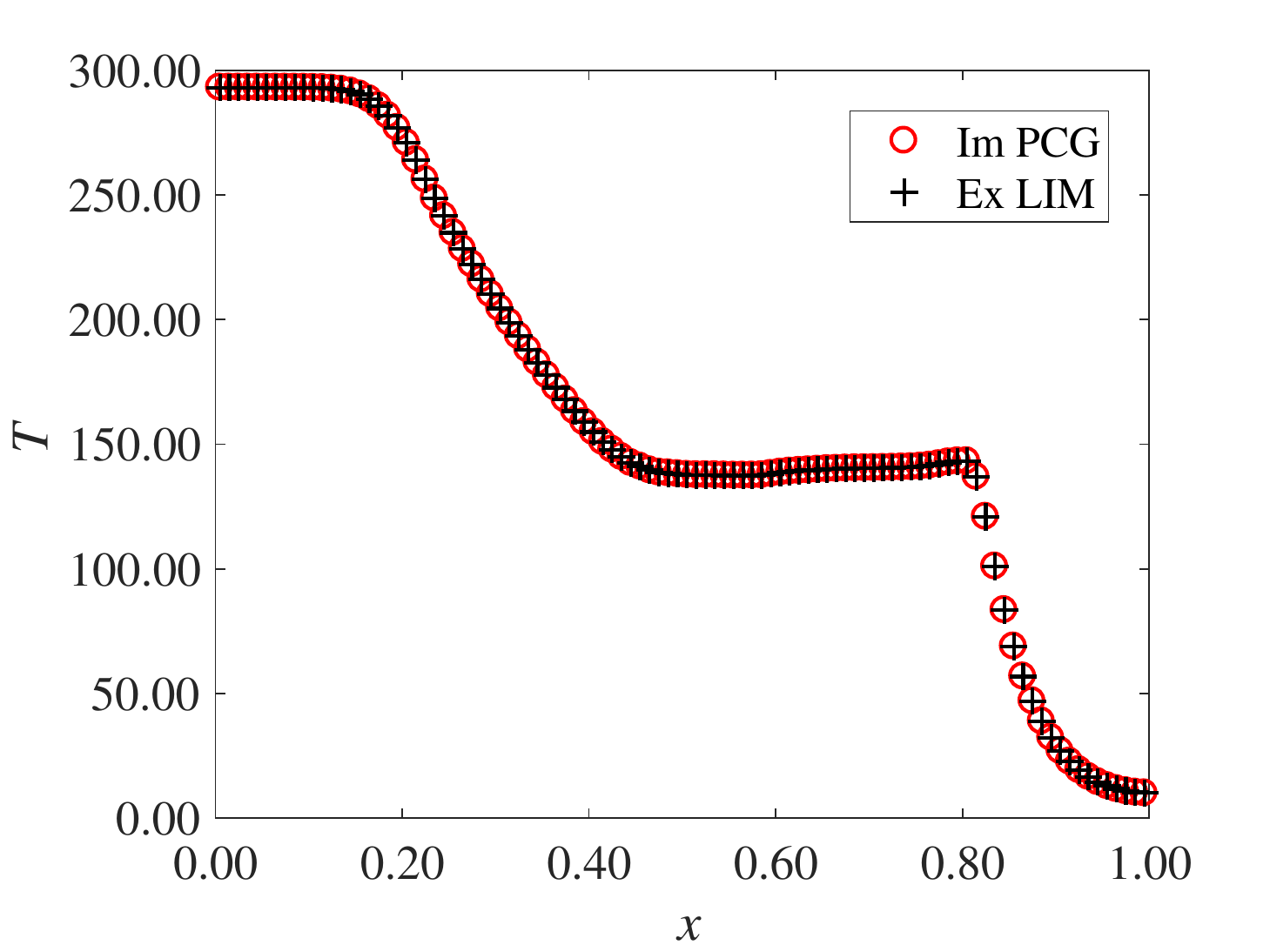}}
\caption{Comparison between the temperature solutions obtained with the implicit PCG method and the explicit LIM method.}
\label{fig:PCG_LIM} 
\end{figure}

In \Cref{fig:PCG_LIM} the temperature solutions (on a grid of 100 computational cells) computed with the implicit PCG method and the explicit LIM are compared. We can see that the solutions almost coincide. However, the explicit LIM demonstrates better parallel efficiency for certain problems, for example,  the 3D supersonic flow in the conduit \cite{zhukov2018}.

\subsection{The laser ablation problem}
\begin{figure}[h!]
\centering
\subfloat[Density]{\label{ABL1D:thickCH01}\includegraphics[width=0.55\textwidth]{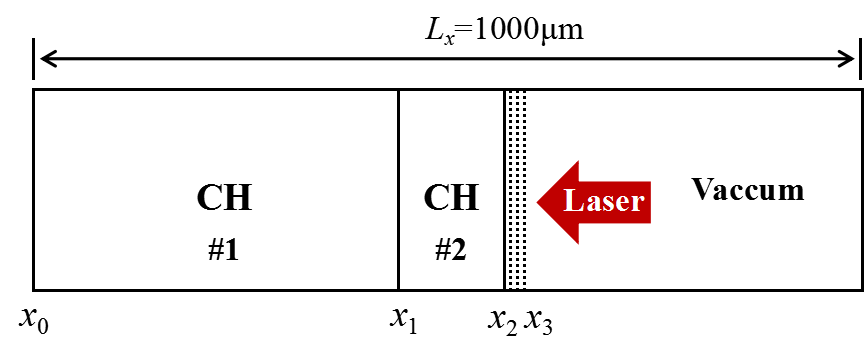}}\\
\subfloat[The absorption area]{\label{ABL1D:thickCH02}\includegraphics[width=0.55\textwidth]{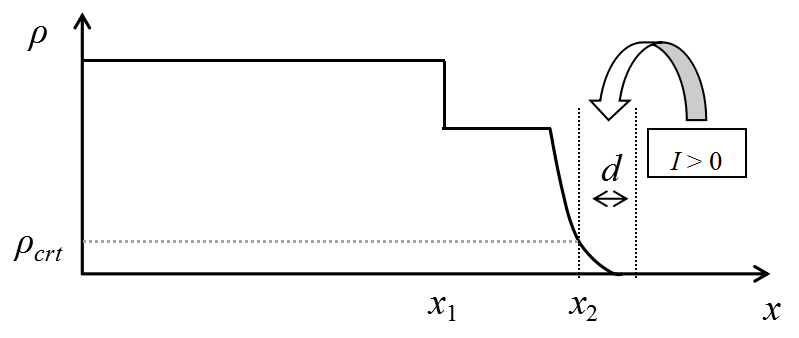}}
\caption{Schematic for the laser ablation problem.}
\label{fig:thickCH0} 
\end{figure}

We consider a problem in the field of ICF -- the laser ablation of a multicomponent target. As displayed in \Cref{ABL1D:thickCH01}, the target is placed in the vacuum. The target consists of two layers of materials which are typically mixtures of CH (phenylethylene ${\text{C}}_8{\text{H}}_8$) and Br (Bromine). Assume that the two layers of CH are characterized by the following parameters:
\[\rho_1 = 1.50 \text{g}/\text{cm}^3, \; \gamma_1 = 2.00, \;\; \text{and} \;\; \rho_2 = 1.00 \text{g}/\text{cm}^3, \; \gamma_2 = 5/3.\]

The specific heats are given to ensure initial pressure and temperature equilibriums,
$C_{v,1} = C_{v,2} = 86.34 \text{cm}^2 / (\mu \text{s} ^2 \cdot \text{MK})$. Extrapolation boundary conditions are imposed on both sides of the computational domain.

The right interface between CH \#2 and the vacuum is planar and is ablated from the right by the laser pulse with wave length  $0.53\mu\text{m}$ and average energy intensity $I = 1.00 \times 10^{14}\text{W}/\text{cm}^2$. The laser energy is absorbed within a distance of $d = 20\mu\text{m}$ to the right of the critical location (\Cref{ABL1D:thickCH02}), where the absorbed energy equals the reflected.  
The critical locations are determined as the rightmost locations where the density is equal to the critical one. The latter is dependent on the chemical properties of of the ablated material and the wavelength of the laser pulse. It is determined as $\rho_{crt} = 1.22\times 10^{-2}\text{g}/\text{cm}^3$ according to the reverse bremsstrahlung theory.  Moreover, we assume that the deposited laser energy intensity is constant within the absorption area $[x_2, x_2 + d]$. The initial density profile is smoothed with an exponential function within $[x_2, x_3]$.
The vacuum on the right is approximated as the material CH \# 2 of very low density, $1.00\times 10^{-5}\text{g}/\text{cm}^3$.

The thermal conductivity is computed with the one-temperature Spitzer-Harm model \cite{Spitzer1953}, 
\begin{equation}
\lambda = 9.44\left(\frac{2}{\pi}\right)^{3 / 2} \frac{\left(k_{B} T_{e}\right)^{5 / 2} k_{B} N_{\mathrm{e}}}{\sqrt{m_{e}} e^{4}} \frac{1}{N_{i} Z_{e}\left(Z_{e}+4\right) \ln \Lambda_{e i}},
\end{equation}
where $k_{B}$ is the Boltzmann constant, $T_{e}$ is the electronic temperature, $N_{e}$ is the electron density, $e$ is the electronic charge, $m_{e}$ is the electronic mass, $N_{i}$ is the ion density, $Z_{e}$ is the atomic number. For a certain plasma,
\begin{equation}
N_{i}=\frac{N_{0}}{A_{c}} \rho, \quad N_{e}=Z_{e} N_{i},
\end{equation}
where $A_{c}$ is the average atomic weight, $N_{0}$ is the Avogadro's number.

$\ln \Lambda_{e i}$ is the Coulomb logarithm of laser absorption and determined with
$$
\ln \Lambda_{e i}=\left\{\begin{array}{ll}
\max \left(1, \ln \frac{l_{D}}{l_{L D}}\right), & \frac{Z_{e}^{2}}{3 k_{B} T_{e}} \geq l_{d B}, \\
\max \left(1, \ln \frac{l_{D}}{l_{d B}}\right), & \frac{Z_{e}^{2}}{3 k_{B} T_{e}}<l_{d B} .
\end{array}\right.
$$
where $l_{D}$ is Debye length, $l_{L D}$ is Landau length, $l_{d B}$ is De Broglie wavelength.

The dynamic viscosity of the plasma is calculated with the formula of Braginskii \cite{Robey2003The}
\begin{equation}
\mu = 3.30 \times 10^{-5} \frac{\sqrt{A_c}{T_e}^{5/2}}{\ln \Lambda_{e i} Z_e^4}.
\end{equation}

The mixture thermal conductivity is averaged with volume fractions, which are absent in the conservative four-equation model. For comparison purpose, we use the the parameters of phenylethylene (i.e., $A_c = 6.5$ and $Z_e = 3.5$)  to compute the thermal conductivity $\lambda_k$ and the dynamic viscosity $\mu_k$  for all components.

\begin{figure}[h!]
\centering
\subfloat[Density]{\label{ABL1D:DENS}\includegraphics[width=0.5\textwidth]{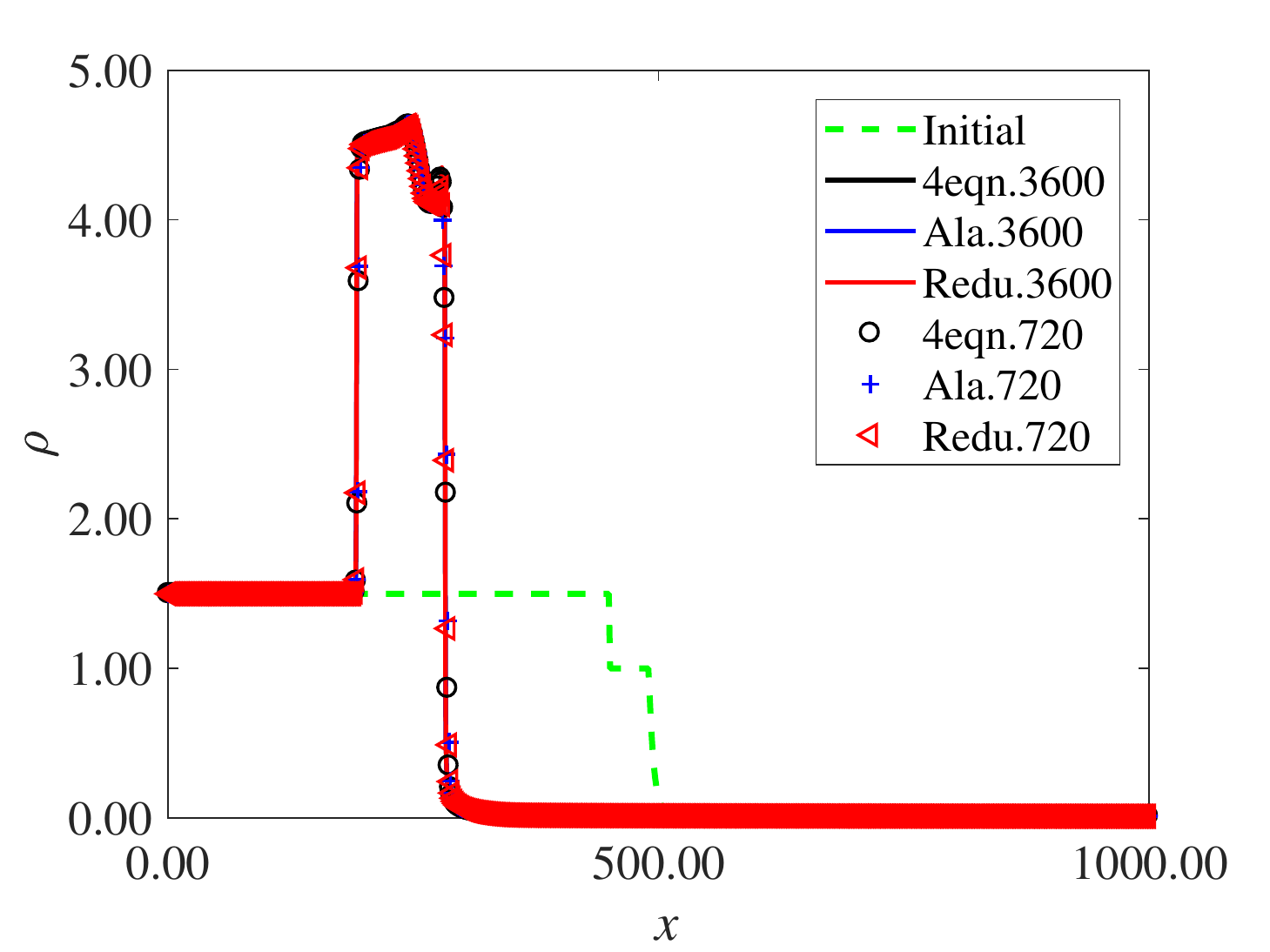}}
\subfloat[Density, locally enlarged]{\label{ABL1D:DENS1}\includegraphics[width=0.5\textwidth]{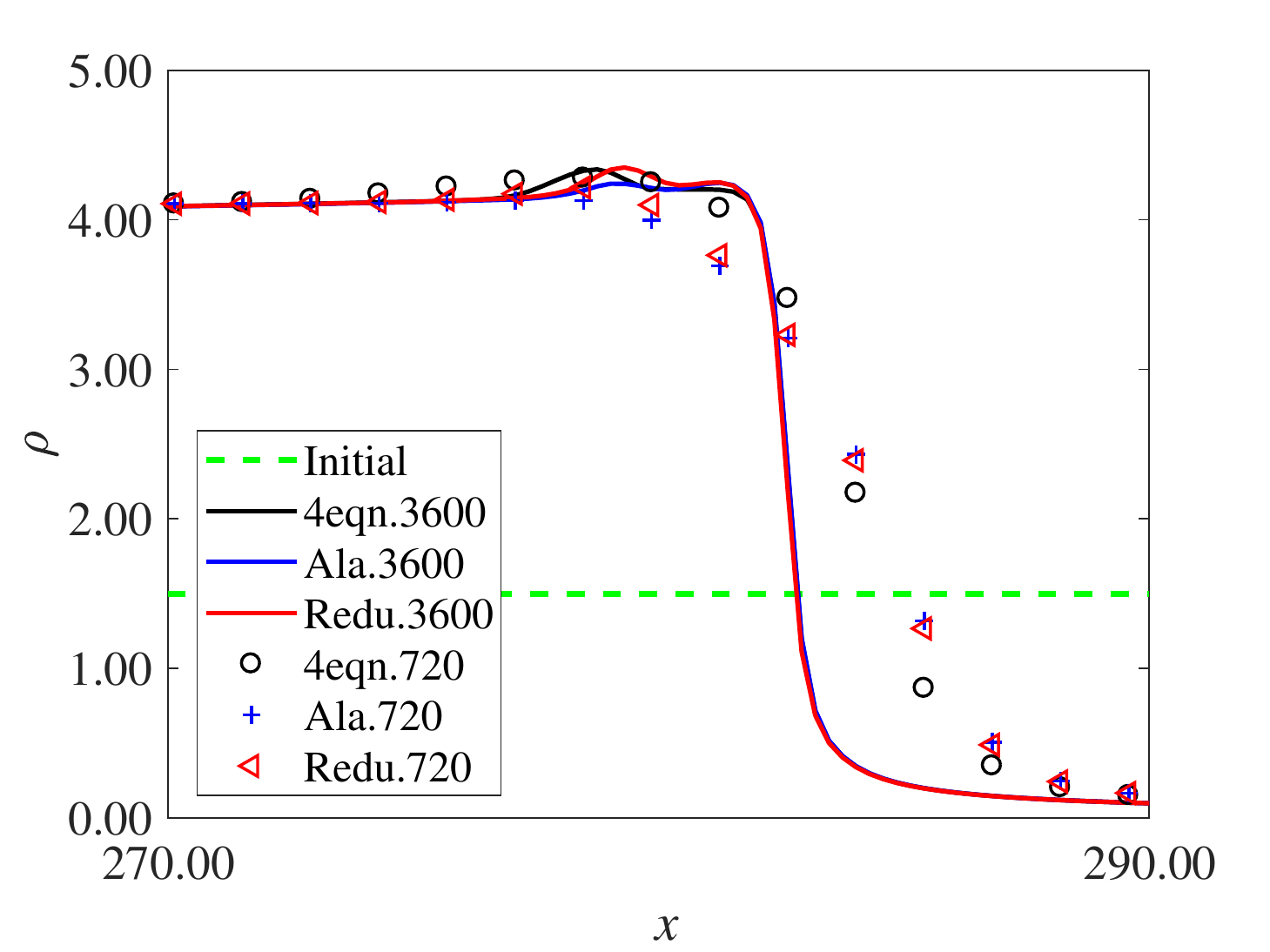}}\\
\subfloat[Pressure]{\label{ABL1D:PRES}\includegraphics[width=0.5\textwidth]{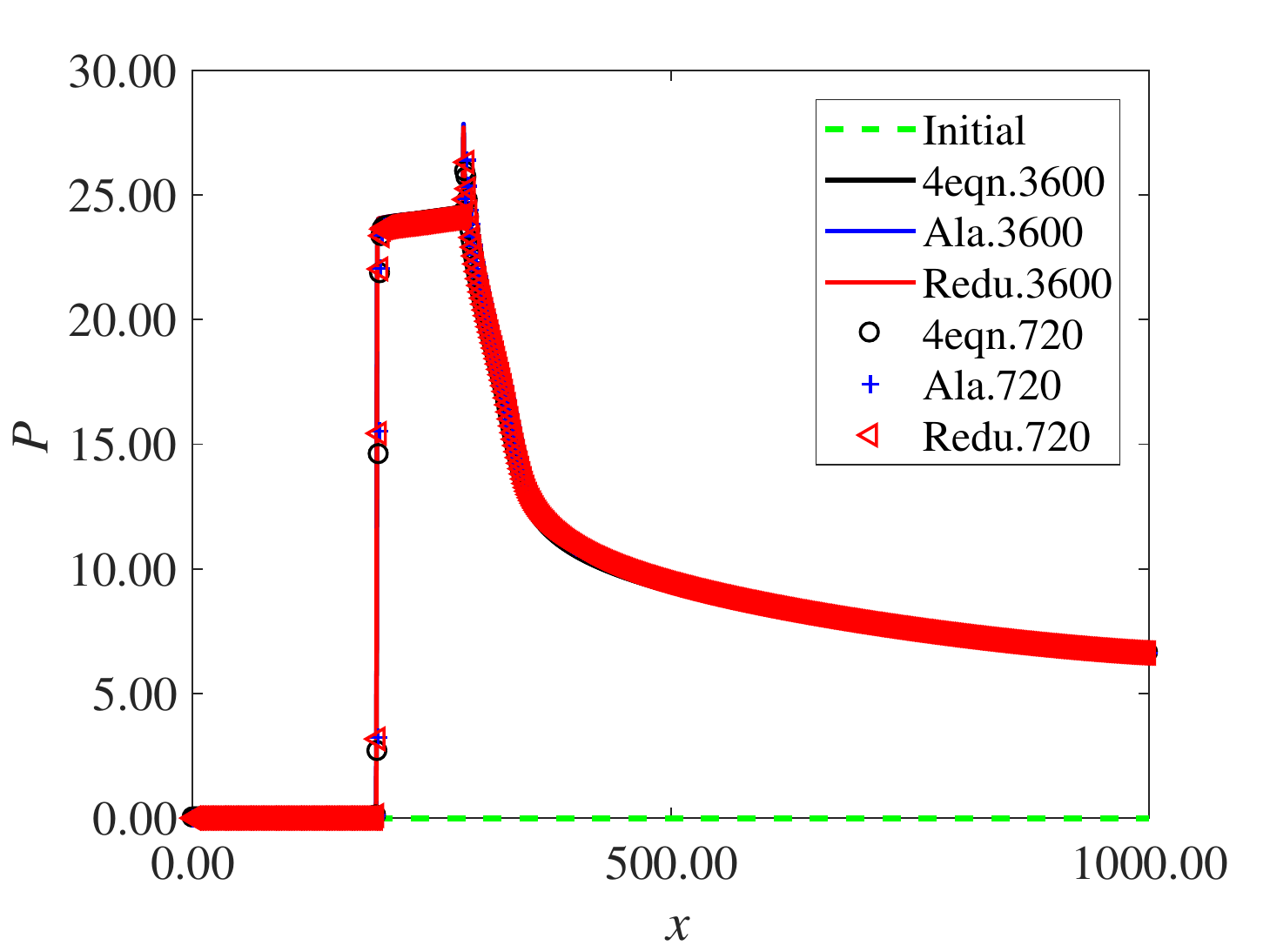}} 
\subfloat[Pressure, locally enlarged]{\label{ABL1D:PRES1}\includegraphics[width=0.5\textwidth]{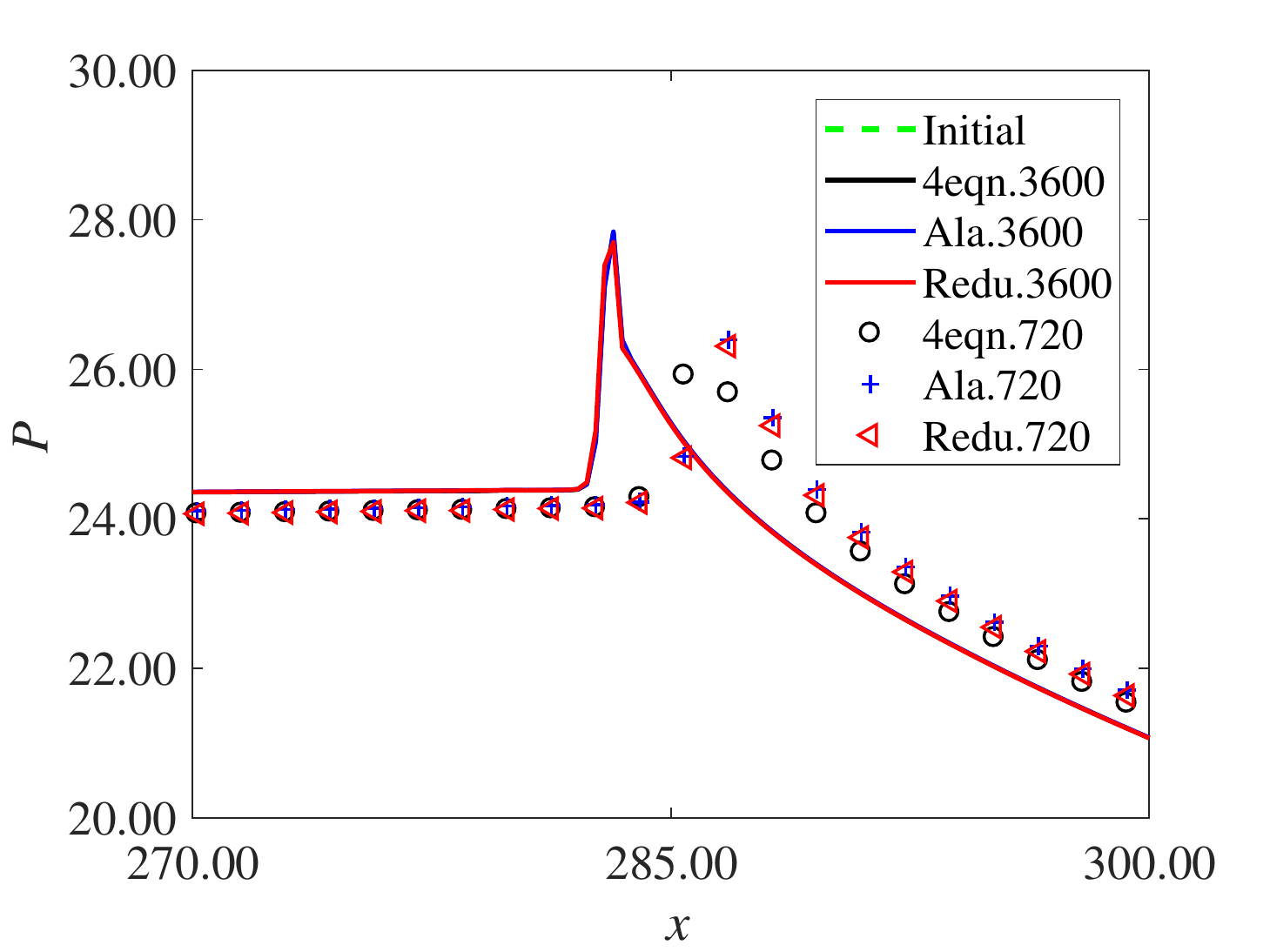}}\\
\subfloat[Velocity]{\label{ABL1D:VEL1}\includegraphics[width=0.5\textwidth]{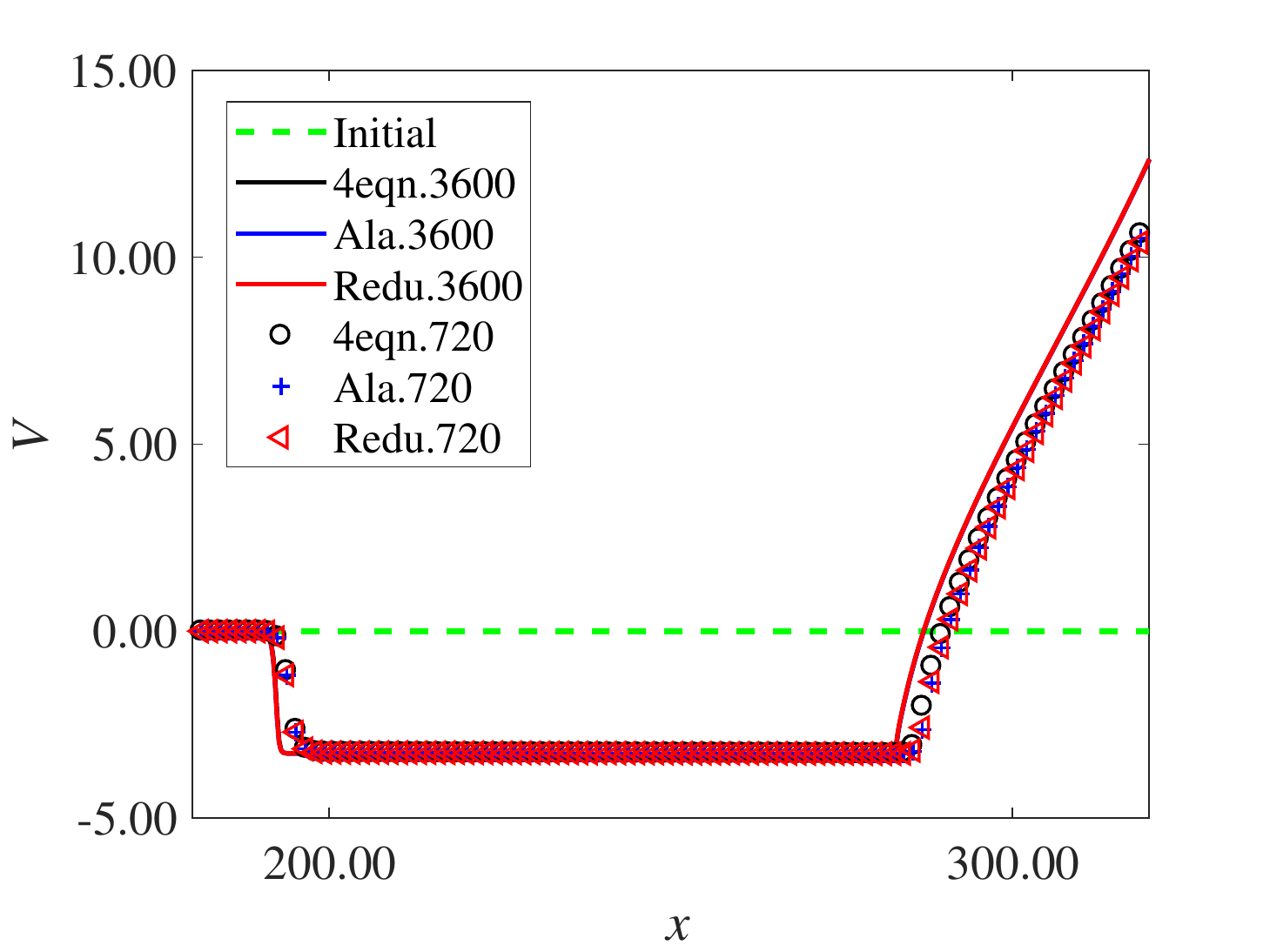}} 
\subfloat[Volume fraction]{\label{ABL1D:VOL1}\includegraphics[width=0.5\textwidth]{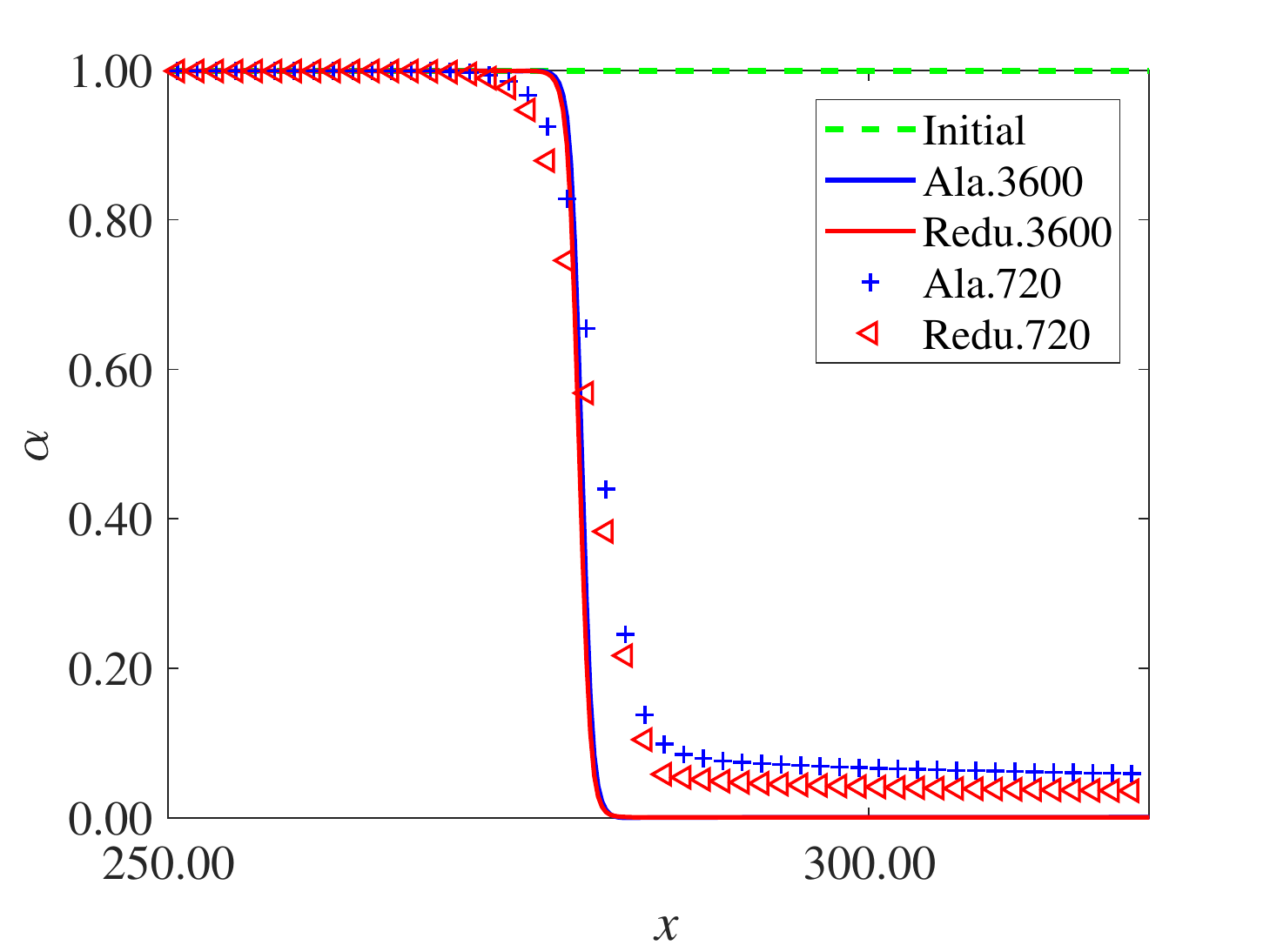}}
\caption{Comparison between the numerical results obtained with different models for the 1D laser ablation problem.}
\label{fig:ABL1D} 
\end{figure}

We perform computations with different models on two uniform grids consisting of 720 and 3600 cells. We choose a larger threshold $\varepsilon = 5.00 \times 10^{-4}$ since with smaller $\varepsilon$ the one-temperature model fails due to spurious oscillations in density. 
The numerical results (for density, pressure, velocity and volume fraction) obtained with different models  are displayed in \Cref{fig:ABL1D}. It can be seen that the numerical results obtained with different models tend to converge to the same solution with the grid refinement. From \Cref{ABL1D:VOL1} one can see that the results for the volume fraction obtained with the one-temperature model suffers from more numerical diffusion than the reduced model.

\begin{figure}[h!]
\centering
\subfloat[Density]{\label{GOOD_KAP:DENS}\includegraphics[width=0.5\textwidth]{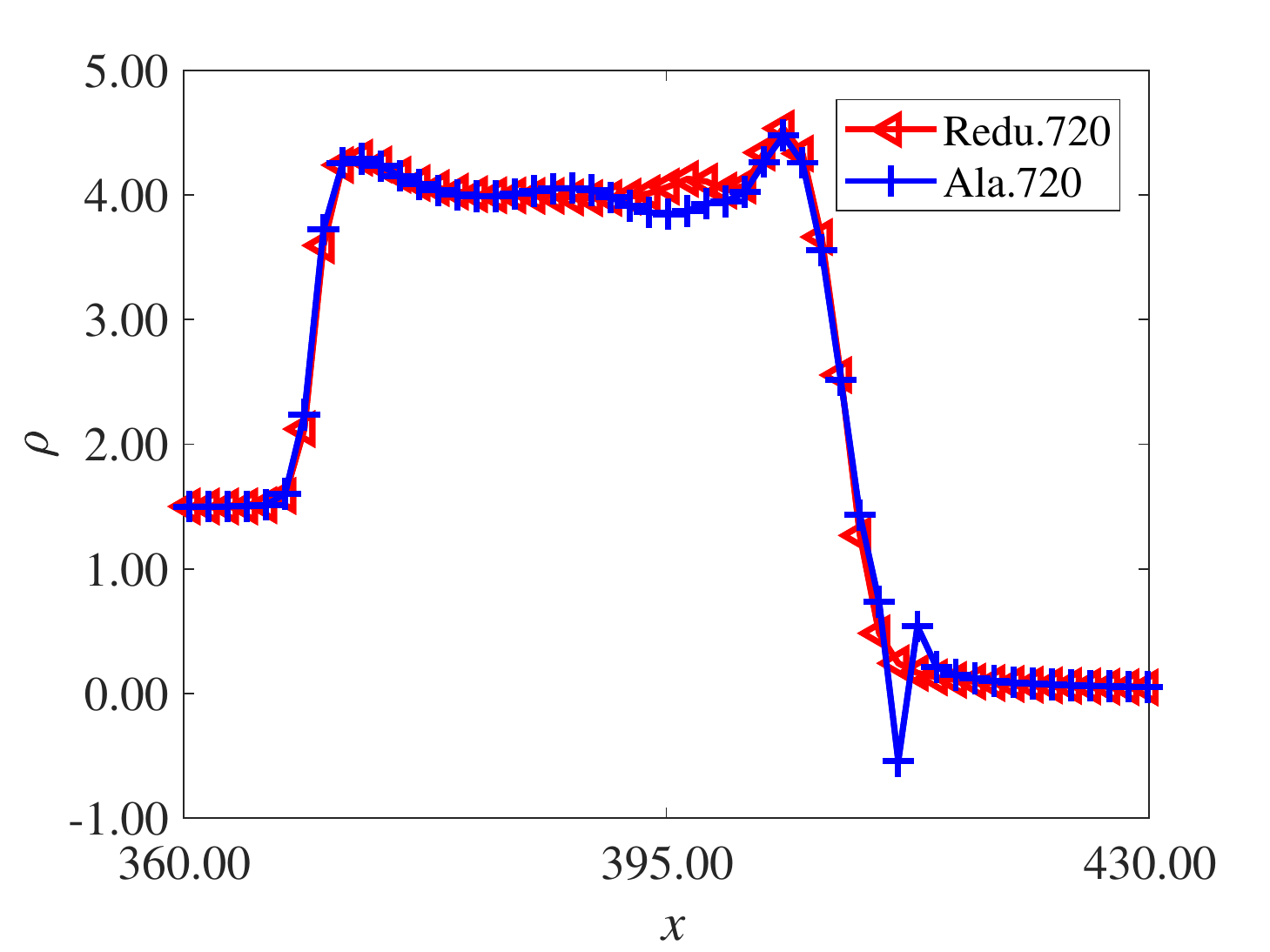}}
\subfloat[Temperature]{\label{GOOD_KAP:DENS1}\includegraphics[width=0.5\textwidth]{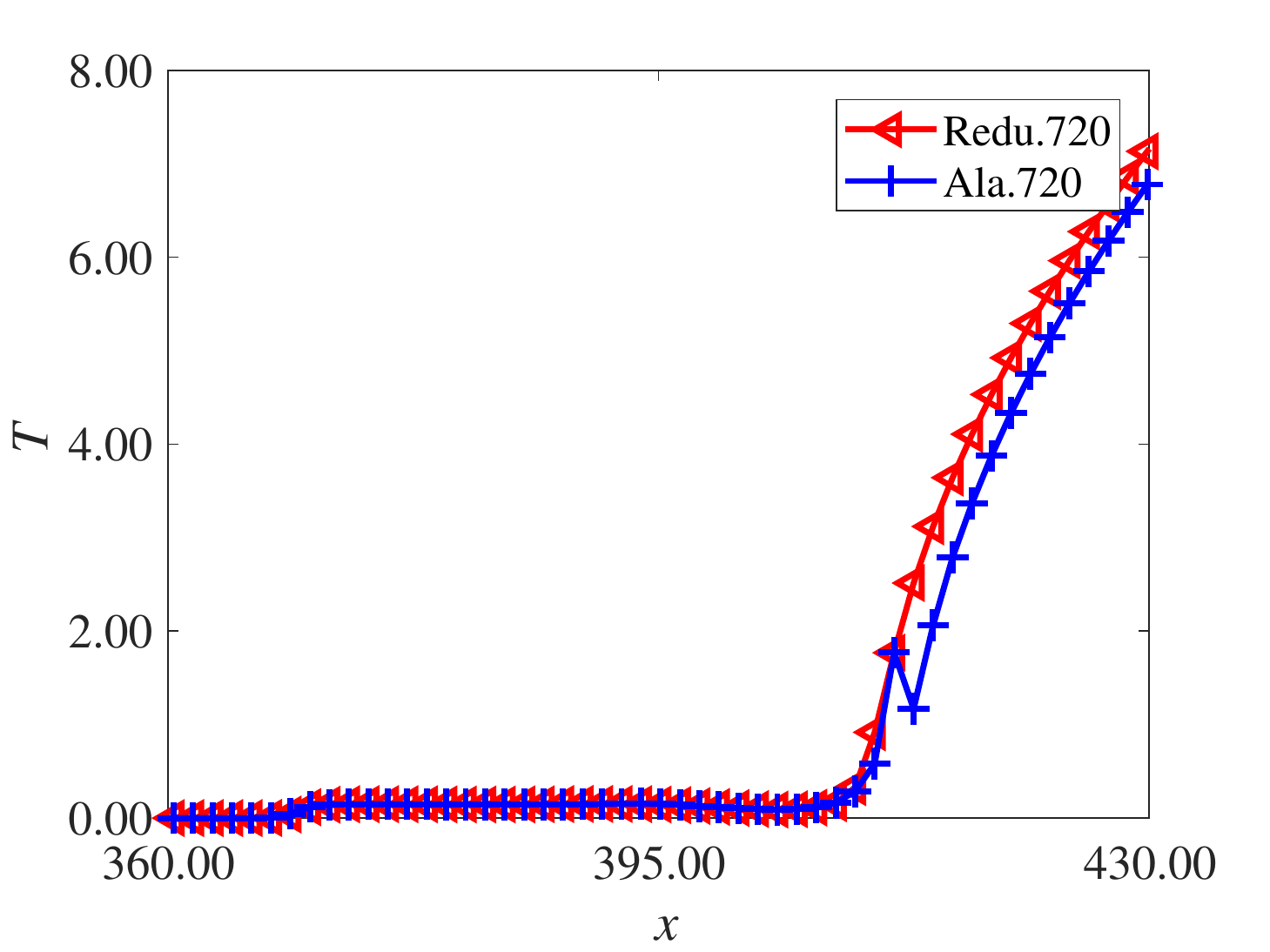}}
\caption{Comparison between the numerical results of different models for the 1D laser ablation problem.}
\label{fig:GOOD_KAP} 
\end{figure}

As mentioned above, our computations show that the one-temperature model is more sensitive to the parameter $\varepsilon$ and yields non-physical solutions with too small $\varepsilon$. The numerical results at $2.49$ns for the density computed by the one-temperature model and the reduced model with $\varepsilon = 1.00 \times 10^{-6}$ are compared in \Cref{fig:GOOD_KAP}. It is clear that the one-temperature model gives non-physical and non-monotone numerical solutions. The emergence of negative density results in the failure of further computation. The reduced model is more robust when $\varepsilon$ becomes smaller and free of this oscillation issue.

\subsection{The triple-point problem}
In this section we consider a benchmark problem for multicomponent flows -- the triple-point problem \cite{Kucharik2010ACS}. Unlike the two-component formulations in literature, we treat this problem as a three-component one. Moreover, we include the diffusion processes. At the initial moment three different states/components are given in different domains as follows: 
\begin{equation}
\left(\rho_{1}, \; \rho_{2}, \; \rho_{3}, \; p, \alpha_{1}, \; \alpha_{2}\right)=\left\{\begin{array}{ll}
\left(1,\;1,\;1,\;1,\;1-\varepsilon, \;\varepsilon/2\right) & \text { if } (x,y) \in \Omega_L, \\
\left(1,\;1,\;1,\;0.1,\;\varepsilon/2, \;1-\varepsilon\right) & \text { if } (x,y) \in \Omega_B, \\
\left(0.125,\;0.125,\;0.125,\;0.1,\; \varepsilon/2, \; \varepsilon/2\right) & \text { if } {{(x,y) \in \Omega_T}. }
\end{array}\right.
\end{equation}
where the left sub-domain $\Omega_L = [0.0,1.0]\times [0.0,3.0]$, the bottom sub-domain $\Omega_B = [1.0,7.0]\times [0.0,1.5]$, the top sub-domain $\Omega_T = [1.0,7.0]\times [1.5,3.0]$. The fluids initially occupy $\Omega_L$, $\Omega_B$ and $\Omega_T$ are referred to as the 1-fluid, the 2-fluid and the 3-fluid.

The three fluids are at rest initially. The adiabatic coefficients for the three ideal gases are $\gamma_1 = 1.5$, $\gamma_2 = 1.4$, and $\gamma_3 = 2.0$, respectively. The  specific heats at constant volume are given as $C_{v,1} = 40.00$,  $C_{v,2} = 50.00$, and $C_{v,3} = 20.00$, which ensures the initial temperature equilibrium. The phase dynamic viscosities are $\mu_1 = 0.10$, $\mu_2 = 0.20$, and  $\mu_3 = 0.05$. The phase coefficients of heat conduction are $\lambda_1 = 0.50$, $\lambda_2 = 1.00$, and $\lambda_3 = 2.00$, respectively.


No-reflecting boundary conditions are imposed on all the boundaries. Computations are performed on a grid of 1400$\times$600 uniform cells. The second-order MUSCL scheme is used for spatial reconstruction, and the Overbee scheme \cite{Chiapolino2017} is used for interface sharpening. 

Direct implementation of the TVD limiters to the volume fractions leads to the violation of the constraint $0 < \alpha_k < 1$ in the case of $N \geq 3$ due to the non-linearity of the limiters.  To overcome this issue, we use the scheme developed in \cite{Zhangchao2020}.

We perform four numerical tests including different physical processes as follows:
\begin{enumerate}
\item[(1)] The hydrodynamic process,
\item[(2)] The hydrodynamic and the viscous processes,
\item[(3)] The hydrodynamic and the temperature relaxation processes, 
\item[(4)] The hydrodynamic, the temperature relaxation, and the heat conduction processes. 
\end{enumerate}

\begin{figure}[h!]
\centering
\label{tripleDens:H}\includegraphics[width=0.95\textwidth]{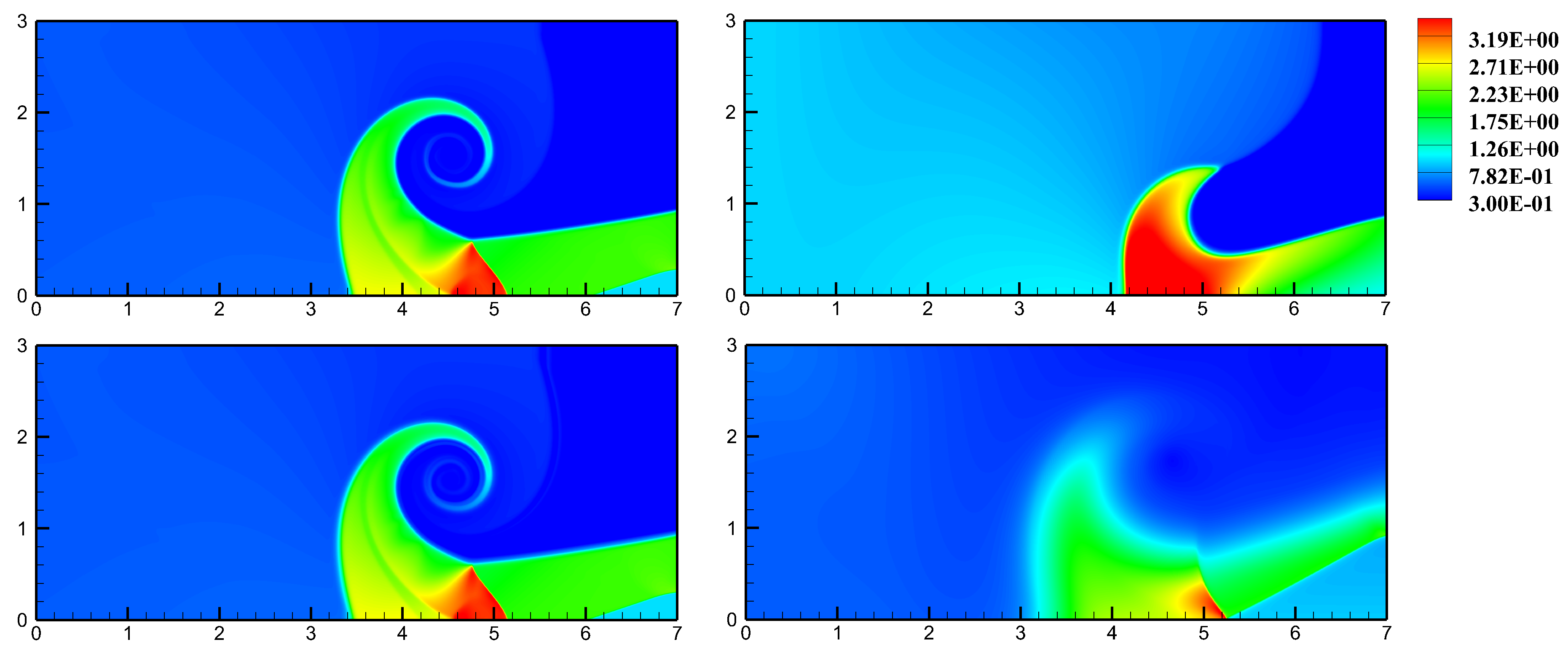}
\caption{Density distributions for the triple point problem with different processes.
Left top: H, right top: H and V, left bottom: H and TR, right bottom: H, TR and HC. Here, the abbreviations used are as follows: 
H -- hydrodynamic, V -- viscous, TR -- temperature relaxation, HC -- heat conduction. The same abbreviations are used for figures in this subsection.}
\label{fig:tripleDens} 
\end{figure}

\begin{figure}[h!]
\centering
\label{tripleSH:H}\includegraphics[width=0.95\textwidth]{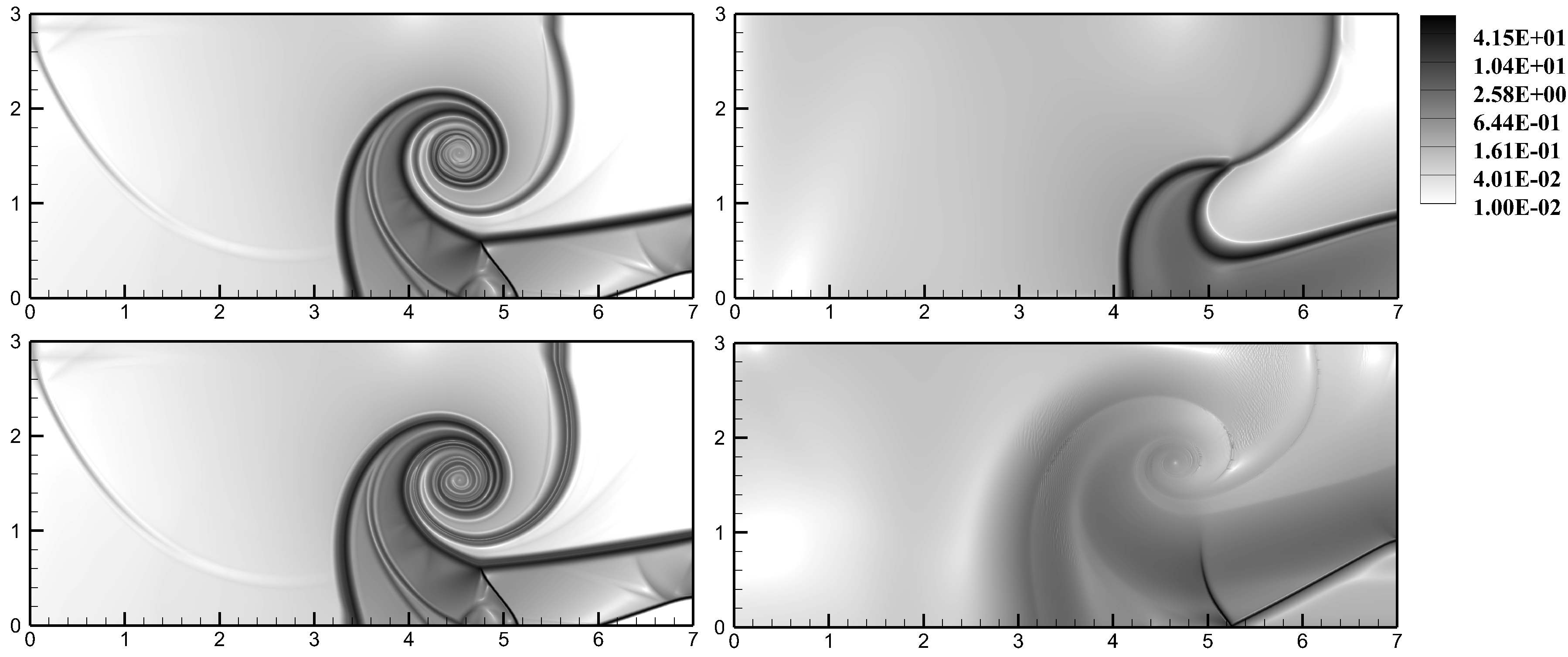}
\caption{Numerical Schlieren for the triple point problem with different processes. Left top: H, right top: H and V, left bottom: H and TR, right bottom: H, TR and HC.}
\label{fig:tripleSH} 
\end{figure}

\begin{figure}[h!]
\centering
\includegraphics[width=0.95\textwidth]{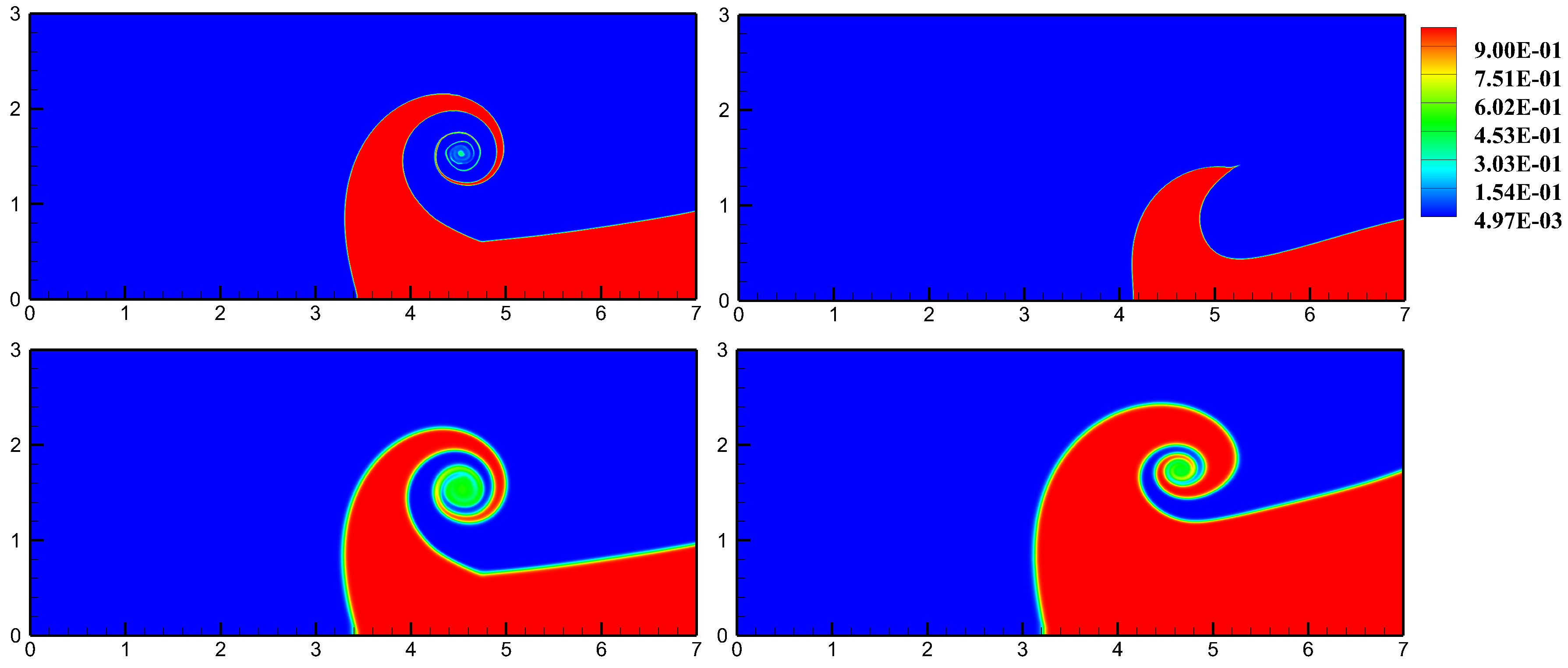}
\caption{Volume fraction of the 2-fluid in the triple point problem with different processes. Left top: H, right top: H and V, left bottom: H and TR, right bottom: H, TR and HC.}
\label{fig:tripleZ2} 
\end{figure}

\begin{figure}[h!]
\centering
\includegraphics[width=0.95\textwidth]{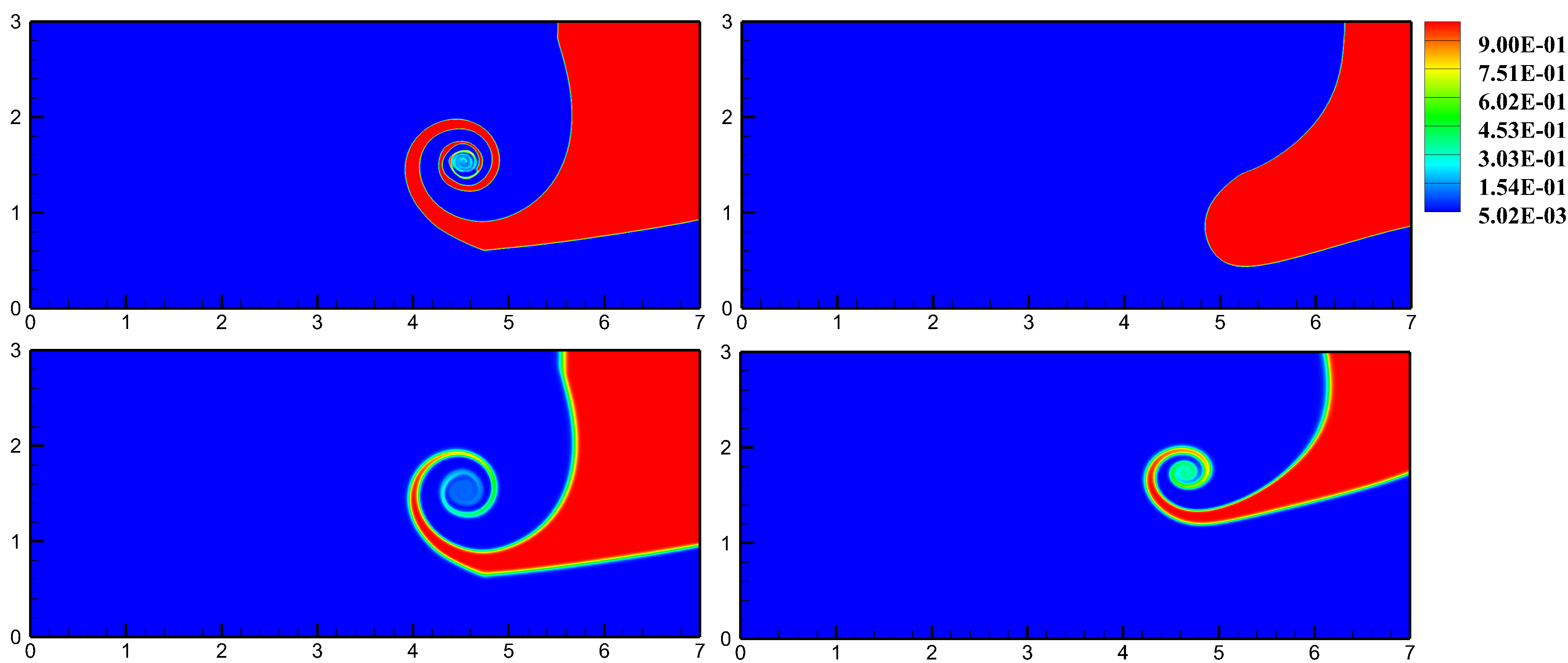}
\caption{Volume fraction of the 3-fluid in the triple point problem with different processes. Left top: H, right top: H and V, left bottom: H and TR, right bottom: H, TR and HC.}
\label{fig:tripleZ3} 
\end{figure}

\begin{figure}[h!]
\centering
\includegraphics[width=0.95\textwidth]{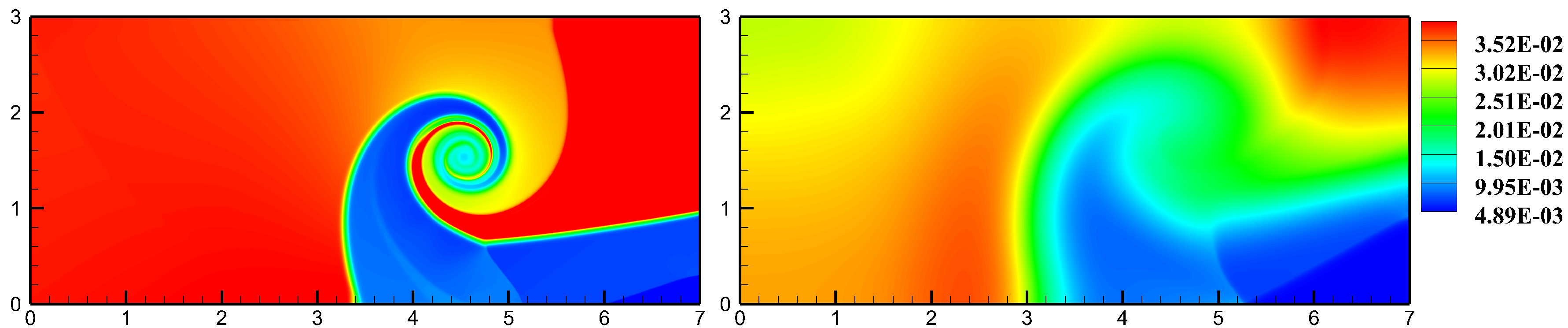} 
\caption{Temperature for the triple point problem with different processes. Left: H and TR, right: H, TR and HC.}
\label{fig:tripleT} 
\end{figure}

The numerical results for density, volume fractions and temperature  of these tests when $t = 5.00$ are displayed in 
\Cref{fig:tripleDens,fig:tripleSH,fig:tripleZ2,fig:tripleZ3,fig:tripleT}.
Shock waves propagate in the 2-fluid and the 3-fluid as a result of the breakup the initial discontinuity between the 1-fluid and the 2-fluid, and the discontinuity between the 1-fluid and the 3-fluid, respectively. Due to the difference in acoustic impedance, the shock waves travel with different velocities. As a result, the Kelvin–Helmholtz instability develops on the interface between the 2-fluid and the 3-fluid.    

Viscosity  diffuses the vorticity and less interface deformation is observed due to the stabilizing effect of the viscosity (\Cref{fig:tripleZ2,fig:tripleZ3}). Moreover, the waves inside the 2-fluid have been significantly smeared (\Cref{fig:tripleSH}) by the viscosity diffusions.

Temperature relaxation drives the phase temperatures into an equilibrium temperature, which is in fact an  temperature averaged in a particular way (see \cref{eq:tr_temp_av}). This temperature averaging procedure results in the smearing of the material interface as shown in \Cref{fig:tripleZ2,fig:tripleZ3}.

As demonstrated in the last sub-figures of \Cref{fig:tripleDens,fig:tripleSH,fig:tripleT}, the heat conduction process significantly smears the density and temperature distribution, however, it barely impacts the thickness of the material interfaces (\Cref{fig:tripleZ2,fig:tripleZ3}).

\subsection{Shock passage through a cylindrical bubble}
\begin{figure}[h!]
\centering
\includegraphics[width=0.6\textwidth]{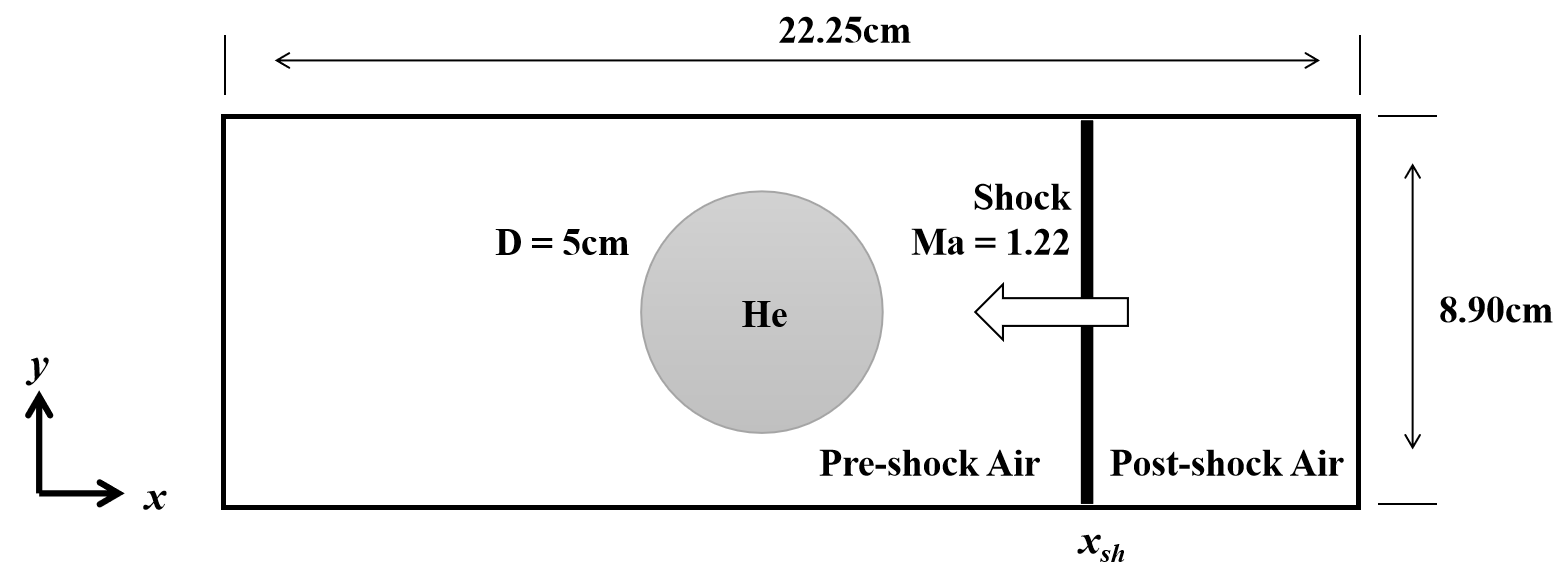} 
\caption{Schematic of the bubble-shock interaction problem.}
\label{fig:HeBubbleSetup} 
\end{figure}

In this section we consider the interaction of shock with a cylindrical helium bubble. The helium inside the bubble is contaminated with 28\% of mass concentration of air. The schematic for this problem is displayed in \Cref{fig:HeBubbleSetup}. The size of the computational domain is 22.25cm $\times$ 8.90cm. The diameter of the bubble is $D = 5.00$cm and the center is (13.80cm, 4.45cm). A left-going  shock wave of Mach 1.22 is initially located at $x_s = 16.80$cm and impacts the helium bubble from the right.

The initial condition $\left(\rho, u, p, \gamma \right)$ are given in the following manner:
\begin{enumerate}
\item[] Post-shock air:  $\left(1.66 \mathrm{~kg} \cdot \mathrm{m}^{-3},-114 \mathrm{~m} \cdot \mathrm{s}^{-1}, 159080.98 \mathrm{~Pa}, 1.4\right)$,
\item[] Pre-shock air: $\left(1.2062 \mathrm{~kg} \cdot \mathrm{m}^{-3}, 0,101325 \mathrm{~Pa}, 1.4\right)$,
\item[] Helium mixture inside the bubble : $\left(0.2204 \mathrm{~kg} \cdot \mathrm{m}^{-3}, 0,101325 \mathrm{~Pa}, 1.6451\right)$.
\end{enumerate}

The specific heats of air and the helium mixture are $717.50$J/(kg$\cdot$K) and $2430.35$J/(kg$\cdot$K), respectively. With the above parameters, the temperatures of air and helium are in equilibrium initially. The pre-shock equilibrium temperature is $T = 293.23$K, and the post-shock one $T = 334.44$K.

We consider this problem with the diffusions being included. The viscosity is calculated with the Sutherland's equation \cite{Sutherland}:
\begin{equation}
\mu = \mu_0 \frac{T_0 + W}{T+W} \left( \frac{T}{T_0} \right)^{\frac{3}{2}},
\end{equation}
where $\mu_0$ is the viscosity at the reference temperature $T_0$. The parameter $W$ is a constant dependent on the material. For air, $\mu_0 = 1.716\times 10^{-5}\text{kg}\cdot \text{m}^{-1} \cdot \text{s}^{-1}$, $T_0 = 273.00$K, $W = 130.00$K; For helium, $\mu_0 = 1.870\times 10^{-5}\text{kg}\cdot \text{m}^{-1} \cdot \text{s}^{-1}$, $T_0 = 273.00$K, $W = 65.00$K. 

The thermal conductivity is determined as 
\begin{equation}
\lambda_{\text{air}} = \frac{\mu_{\text{air}} C_{p,\text{air}}}{Pr},
\end{equation}
Here, as in \cite{2018Capuano}, we use $Pr = 0.7$.

The thermal conductivity of helium is calculated with the following polynomial fitting of experimental data \cite{2018Capuano}:
\begin{equation}
\lambda_{\text{He}} = 1.2900 \times 10^{-11}T^3 -7.4500 \times 10^{-8}T^2 + 3.8960 \times 10^{-4}T + 3.7220 \times 10^{-2}.
\end{equation}

\begin{figure}[h!]
\centering
\includegraphics[width=0.95\textwidth]{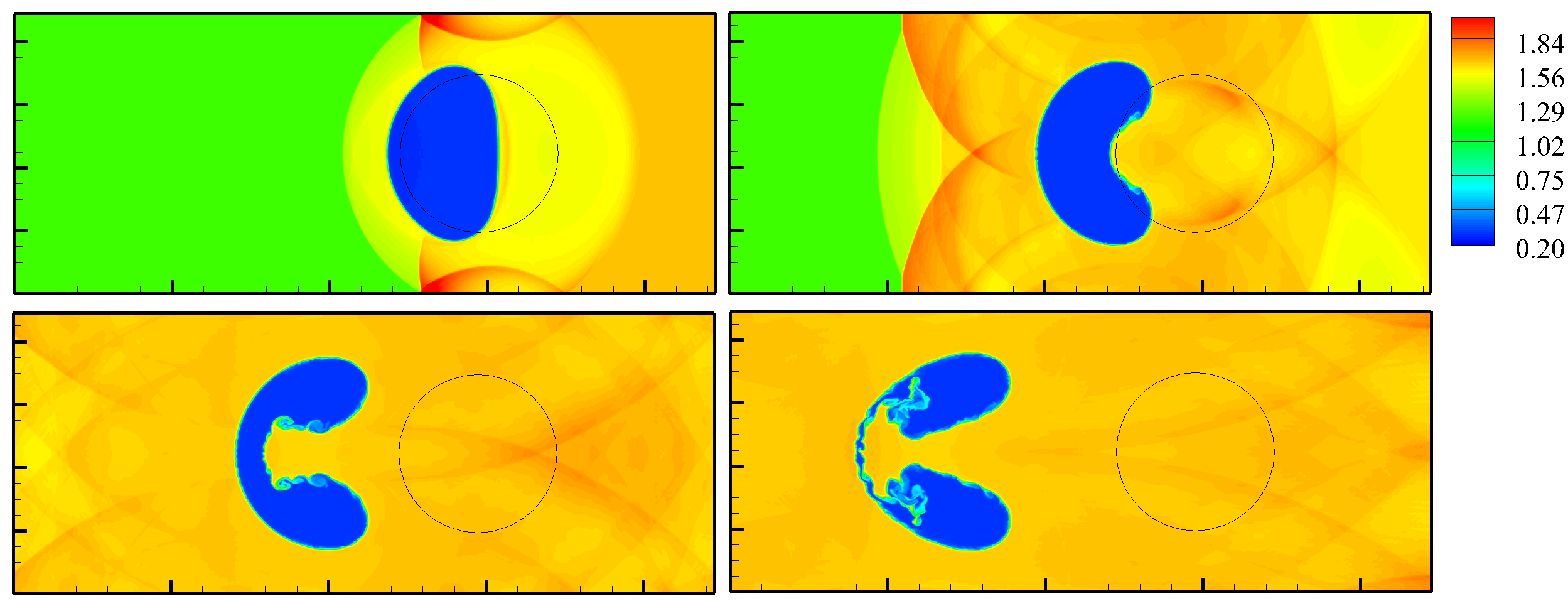} 
\caption{Numerical results for density in the shock-bubble interaction problem when $t = 102, \; 245, \; 427, 674 \mu$s. The time moment $t = 0$ corresponds to the moment when the first contact between the incident shock and the bubble interface happens. The dotted circle represents the initial position of the bubble.}
\label{fig:HeDens} 
\end{figure}

\begin{figure}[h!]
\centering
\includegraphics[width=0.95\textwidth]{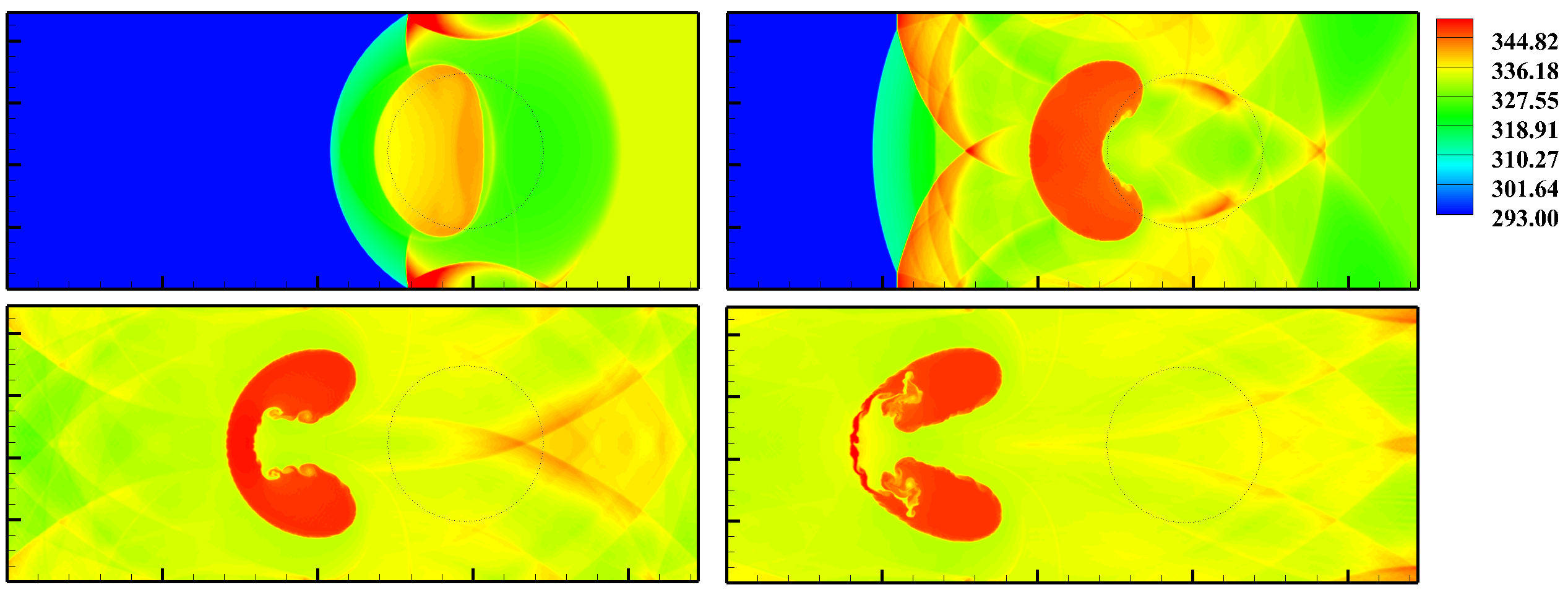} 
\caption{Numerical results for temperature in the shock-bubble interaction problem.}
\label{fig:HeTemp} 
\end{figure}

\begin{figure}[h!]
\centering
\includegraphics[width=0.95\textwidth]{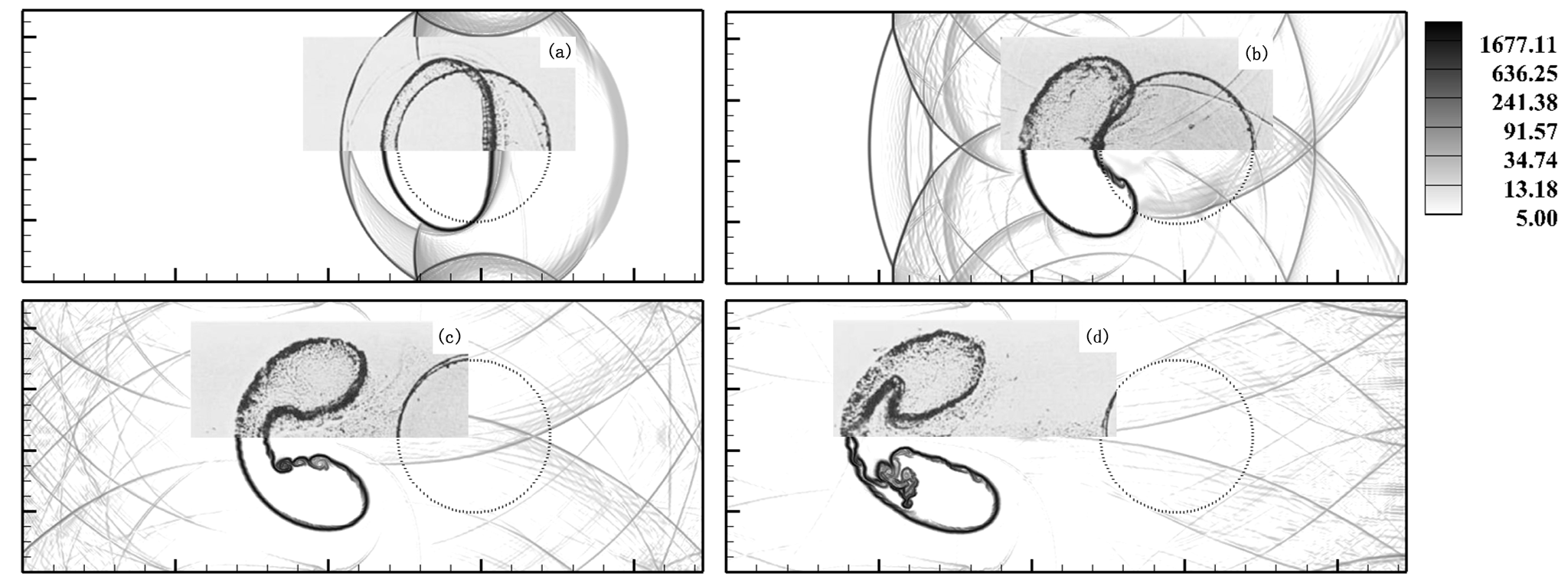} 
\caption{Evolution of the numerical Schlieren with time in the shock-bubble interaction problem. The shaded parts on the top are experimental shadow-photographs from \cite{Haas1987}}
\label{fig:HeSH} 
\end{figure}

Extrapolation boundary conditions are imposed on the left and right boundaries, and periodical conditions on the top and bottom. We perform computation on a uniform mesh of 1200$\times$480 computational cells. The fifth order WENO scheme \cite{Coralic2014Finite,Shu1998Essentially} is used for spatial reconstruction. The numerical results for density, temperature and numerical Schlieren are displayed in \Cref{fig:HeDens}, \Cref{fig:HeTemp}, and  \Cref{fig:HeSH}, respectively.

The bubble undergoes deformation during its interaction with the incident shock, which is reflected and refracted in this process. On the bubble interface Richtmyer–Meshkov (RM) instability develops. For detailed analysis of this process, one may refer to the previous works reviewed in \cite{Ranjan2011}. From \Cref{fig:HeTemp}, one can see that the temperature of the helium mixture inside the bubble is increasing in the course of compression. The numerical Schlieren is compared with the experimental results that are available from \cite{Haas1987}. It can be seen that our numerical results capture the basic deformation behaviours of the bubble.

\subsection{The 2D laser ablation Rayleigh–Taylor instability problem}
\begin{figure}[h!]
\centering
\subfloat{\includegraphics[width=0.45\textwidth]{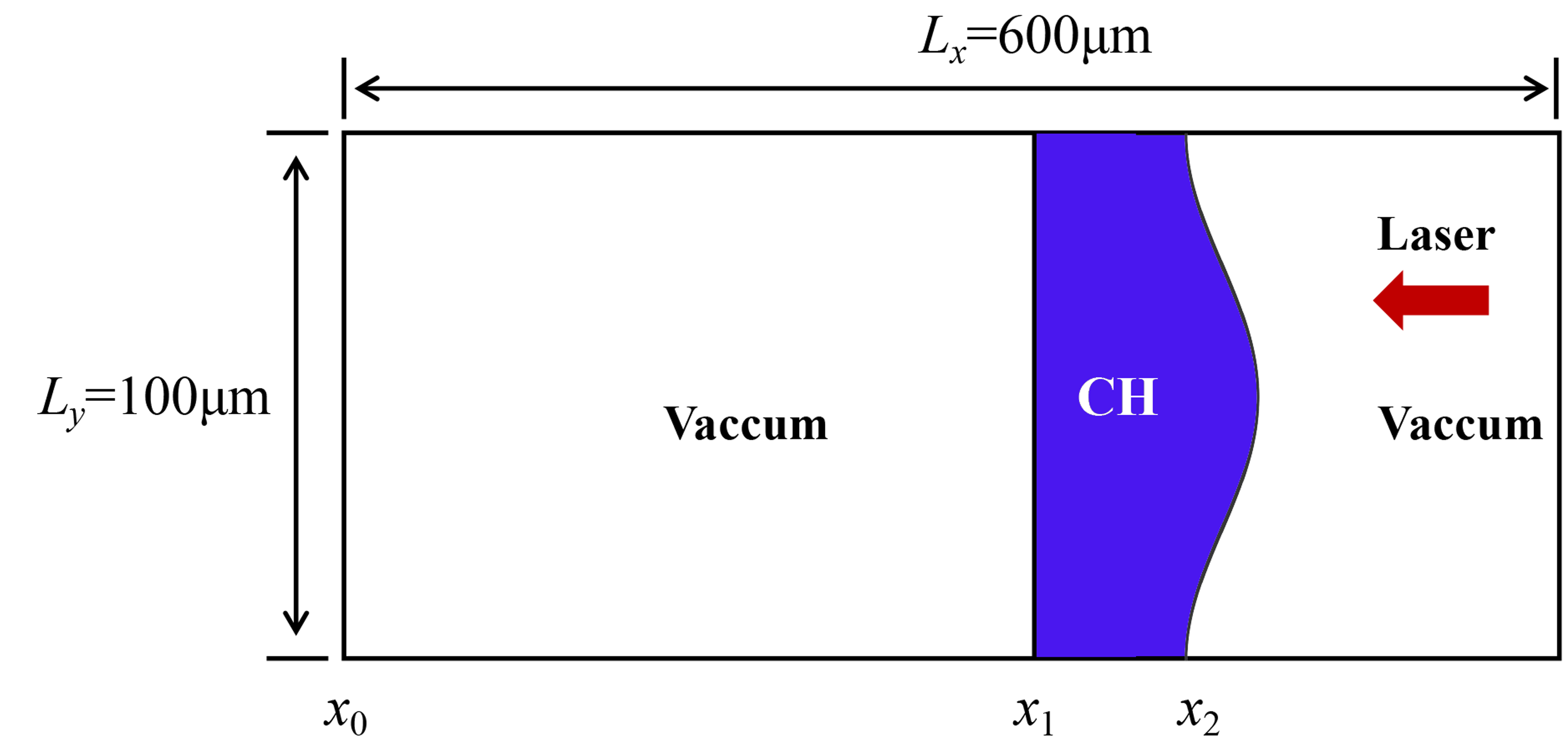}}
\subfloat{\includegraphics[width=0.45\textwidth]{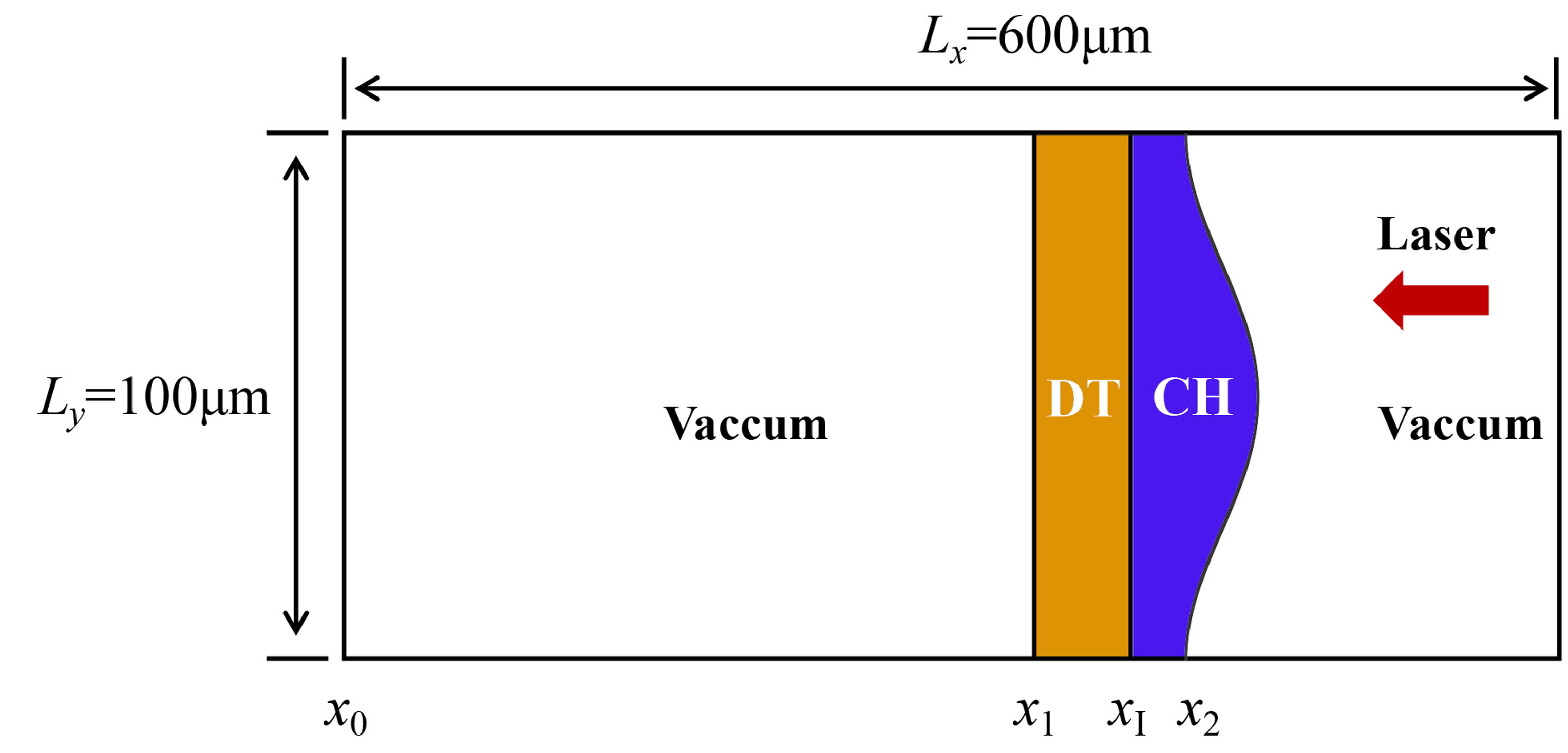}}
\caption{Problem setup of the 2D laser acceleration of corrugated target. Left: one-component target, right: two-component target. The initial positions of the target: $x_1 = 450\mu{\text{m}}, \; x_I = 460\mu{\text{m}},$ and the right interface is perturbed as $x_2 = 475 - 1.5 \text{cos}\left( 2\pi y/L_y \right)$.}
\label{fig:ART_SETUPS} 
\end{figure}

To verify our code, we first consider a one-component problem, then continue with a two-component one.
The set-up of the one-component problem is similar to that in \cite{Li2020Numerical,Nakai1995,Atzeni2004} and illustrated on the left of \Cref{fig:ART_SETUPS}. The  CH foil target is of thickness $25\mu{\text{m}}$.  The rightmost surface of the target is sinusoidally modulated to be a cosine curve with the wavelength $100\mu{\text{m}}$ and the corrugation amplitude $1.5\mu{\text{m}}$. The laser pulse (of wavelength $0.53\mu{\text{m}}$ and energy intensity $2\times 10^{14}\text{W}/\text{cm}^2$) ablates the corrugated interface from the right, resulting in the target acceleration. The initial perturbations are amplified and Rayleigh–Taylor instability (RTI) develops in the vicinity of the right interface. 

We perform simulations on a mesh of 960$\times$160 cells. Periodical boundary conditions are imposed on the horizontal boundaries, and extrapolation conditions on the vertical ones.  Our simulation results for the one-component problem are qualitatively compared against the experimental ones from  \cite{Nakai1995}. It can be seen that the numerical results well capture the basic physics of this process, including the bubble/spike growth and the shell breaking.

  We then consider a two-component target consisting of two layers of different components: CH and DT (Deuterium-Tritium).  The thermodynamic properties of the materials are 
  \begin{enumerate}
  \item[CH ] $\gamma = 1.67, \;\; \rho = 1.00{\text{g}/\text{cm}^3}, \;\; C_v = 86.34{\text{cm}^2 / \left( {\mu\text{s}^2} \cdot \text{MK} \right)}$,
  \item[DT ] $\gamma = 1.41, \;\; \rho = 0.50{\text{g}/\text{cm}^3},\;\;  C_v = 134.70{\text{cm}^2 / \left( {\mu\text{s}^2} \cdot \text{MK} \right)}$.
  \end{enumerate}
  The schematic of the problem set-up is demonstrated on the right of  \Cref{fig:ART_SETUPS}.
  
  The numerical results of the one-component problem and the two-component one (obtained at $2.20$ns and $2.70$ns) are compared in \Cref{fig:PURE_2COM}. It can be seen that the bubble and spike grow faster than those in the one-component problem. In comparison with the pure CH target, the less-weighted composite target achieves larger acceleration with the same amount of absorbed energy. The material interface between CH and DT are well resolved in the course of the computation.
  
\begin{figure}[h!]
\centering
{\includegraphics[width=0.70\textwidth]{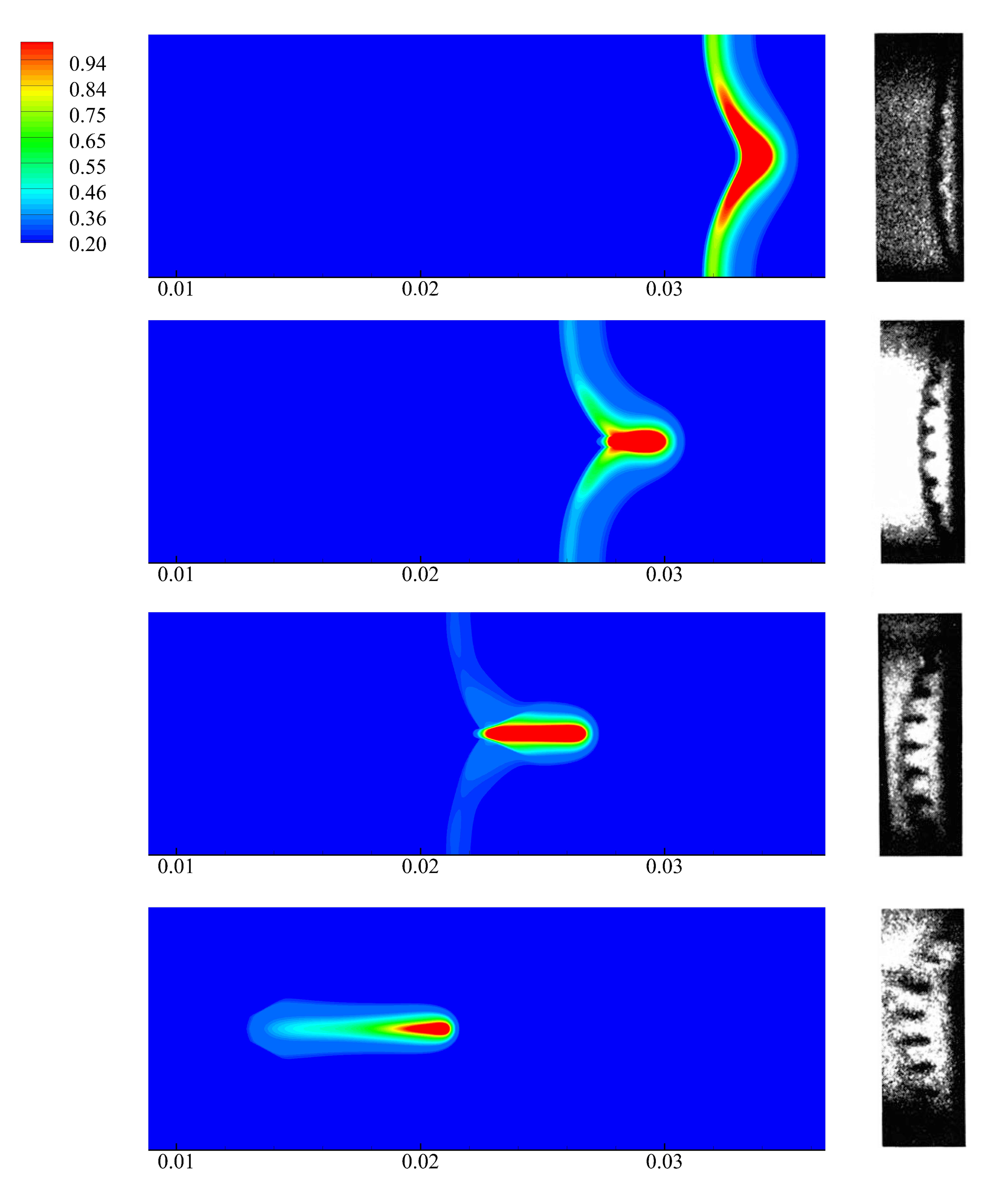}}
\caption{The temporal evolution of the density field in the one-component laser ablation RTI problem. The figures on the right are experiment results  from \cite{Atzeni2004}, where dark regions correspond to the CH target.}
\label{fig:ART2D} 
\end{figure}

\begin{figure}[h!]
\centering
{\includegraphics[width=0.70\textwidth]{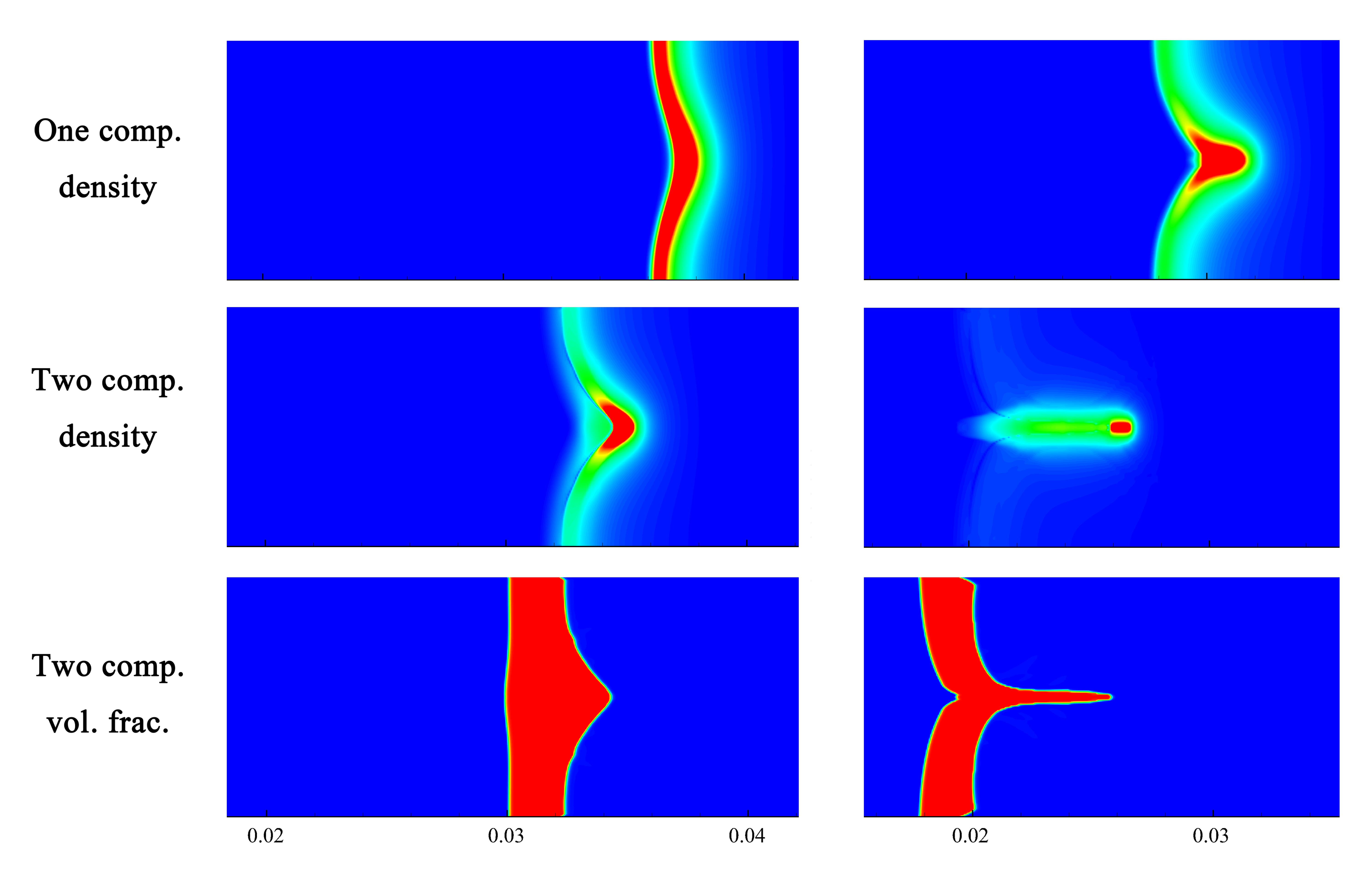}}
\caption{Simulation results for the ablation of the one- and two-component target at 2.2ns (left column) and 2.7ns (right column). The first and second row correspond to the density field of the one- and two-component problem, respectively. The last row presents the DT volume fraction.}
\label{fig:PURE_2COM} 
\end{figure}

\section{Conclusion}
\label{sec:conclusion}
In the present paper we have presented a temperature non-equilibrium model with inter-phase heat transfer, external energy source, and diffusion. We give three equivalent formulations (including the one energy equation formulation and the multiple energy equation formulation) of this model that are reduced from the BN type model for $N$-phase flows in the limit of instantaneous mechanical relaxations. The proposed model ensures the correct interface jump conditions without introducing spurious oscillations in pressure and temperature. Numerical methods for its solution have been proposed on the basis of the fractional step method. The model is split into four sub-systems including the hydrodynamic part, the viscous part, the temperature relaxation part, and the heat conduction part. For solving the hydrodynamic part, the one energy equation is used.  The multiple energy equation formulation is more straightforward to consider inter-phase energy exchange and is thus used for solving the rest parts. The hyperbolic equations involved are solved with the Godunov finite volume method, and the parabolic ones with the locally iterative method based on Chebyshev parameters. 

The developed model and numerical methods have been used for solving several multicomponent problems and verified against analytical and experimental results. Comparisons with results of the one-temperature models show that our model has advantage in convergence and robustness thanks to its physical consistency. Especially, we also have shown the ability of our model in simulating the laser ablation process of a multicomponent target in the ICF field, where the heat conduction plays a significant role.



\normalem
\bibliographystyle{plain} 
\bibliography{ref}



\end{document}